\documentclass{amsart}
\usepackage{amssymb,amsfonts}

\theoremstyle{plain}
\newtheorem{theorem}{Theorem}
\newtheorem{corollary}{Corollary}

\newtheorem{proposition}{Proposition}

\theoremstyle{definition}
\newtheorem{definition}{Definition}

\theoremstyle{remark}
\newtheorem*{notation}{Notation}
\newtheorem*{example}{Example}

\numberwithin{equation}{section}
\newcommand{\GFE}{\mathbb{GFE}}
\newcommand{\flg}{\mathcal{F}(\mathcal{S})[G]}
\newcommand{\fl}{\mathcal{F}(\mathcal{S})}
\newcommand{\f}{\mathbb{F}}

\newcommand{\re}{\mathbb{R}}
\newcommand{\fc}{\mathcal{F}}
\newcommand{\set}{\mathcal{S}}
\newcommand{\gr}{\mathcal{S}}
\newcommand{\e}{\mathcal{E}}
\newcommand{\s}{\mathfrak{S}}
\newcommand{\imp}{\Rightarrow}
\newcommand{\tpk}{\otimes_{k}}
\newcommand{\tpc}{\otimes_{\mathbb F}}
\newcommand{\ds}{\oplus}
\newcommand{\sub}{\subset}
\newcommand{\un}{\bigcup}
\newcommand{\sn}{\bigcap}
\newcommand{\gm}{\varGamma}
\newcommand{\ft}{\longrightarrow}
\newcommand{\mt}{\longmapsto}
\newcommand{\iso}{\approx}
\newcommand{\zc}{\mathbb{C}}
\newcommand{\df}{Difn_*}

\newcommand{\be}{\begin{equation}}
\newcommand{\ee}{\end{equation}}
\newcommand{\bes}{\begin{equation*}}
\newcommand{\ees}{\end{equation*}}
\newcommand{\bpr}{\begin{proposition}}
\newcommand{\epr}{\end{proposition}}
\newcommand{\bp}{\begin{proof}}
\newcommand{\ep}{\end{proof}}
\newcommand{\bas}{\begin{align*}}
\newcommand{\eas}{\end{align*}}

\begin{document}
\title[The theory of linear G-difference equations ]
      {The Theory of linear G-difference equations }
\author{Per K. Jakobsen\\
	Valentin V. Lychagin}
\address{Faculty of science \\
         University of Troms\o\/\\
         Troms\o\/, 9037} 
\email{perj@math.uit.no}
\keywords{Finite difference equations,modules,morphisms,categories}
\subjclass{Primary: 39A05; Secondary: 39A70}
\date{December 17l, 1997}
\begin{abstract}
We introduce the notion of difference equation defined on a structured set. The symmetry group
of the structure determines the set of difference operators. All main notions in the
theory of difference equations are introduced as invariants of the action of the symmetry group.
Linear equations are modules over the skew group algebra, solutions are morphisms relating a given
equation to other equations, symmetries of an equation are module endomorphisms and conserved structures
are invariants in the tensor algebra of the given equation. 

We show that the equations and their solutions
can be described through representations of the isotropy group of the symmetry group of the underlying set.
We relate our notion of difference equation and solutions to systems of classical difference equations and their
solutions and show that out notions incluse these as a special case.
\end{abstract}

\maketitle

\newpage

\tableofcontents

\section{Introduction} \label{I}

Let us consider a general second order difference equation of the form

\begin{align*}
a_if_{i+1}+b_if_i+c_if_{i-1}=0
\end{align*}

Introduce the simple graph $\gr$ consisting of vertices $\{x_i\}_{i\in \mathbb{Z}}$ and edges
$\{\{x_i,x_{i+1}\}\}_{i \in \mathbb{Z}}$. Let $\fl$ be the $\re$-algebra of $\re$-valued functions on
the graph $\gr$. Then the sequences $\{a_i\},\{b_i\},\{c_i\}$ and $\{f_i\}$ are all elements in $\fl$.
Denote these elements by a,b,c,and f. Let s be the operator of left translation on the lattice $\gr$,
$sx_i=x_{i-1}$.  Then s acts on $\fl$ in the natural way
\begin{align*}
(sf)(x_i)=f(s^{-1}x_i)
\end{align*}

Define $\triangle=as+be+cs^{-1}$ where $e$ acts as the  identity on $\gr$. Then $\triangle$ acts on $\fl$ 
as a $\re$-linear operator and our original equation can be written
\begin{align*}
\triangle(f)=0
\end{align*}
In order to understand what $\triangle$ is in algebraic terms we need to introduce some new notions.
Let $G=Aut(\gr)$ be the automorphism group of the graph $\gr$. This group acts on $\fl$ in the natural way
$(gf)(x_i)=f(g^{-1}x_i)$. Let A be the set of finite formal linear combinations of elements in $G$
with coefficients in $\fl$.
\begin{align*}
A=\{\sum_ga_gg\mid a_g\in\fl\}
\end{align*}

On the set A we define addition and scalar multiplication with elements $r\in\re$ componentwise. Product is
 defined in the following way
\begin{align*}
(ag)(bg')=(agb)(gg')
\end{align*}

With these operations A is a $\re$-algebra. $A=\fl[G]$ is the skew group algebra of $G$ over $\fl$. This algebra
acts on $\fl$ through
\begin{align*}
(\sum_ga_gg)f=\sum_ga_gg(f)
\end{align*}

Using these notions we observe that our classical difference operator $\triangle=as+be+cs^{-1}$ is 
an element of the skew group algebra A. It is now evident that we can interpret all elements in
A as difference operators over $\gr$. We will in fact define A to be the algebra of $G$-difference operators
over $\gr$. This means that the notion of difference operator is defined in terms of the symmetries
of the underlying graph $\gr$. The group of symmetries of $\gr$ measures the arbitrariness in the
description of $\gr$. Without this arbitrariness difference operators  could not
exist, in a totally asymmetrical space with trivial symmetry group there could be no difference
operators and as a consequence no difference equations.

We will in this paper generalize these simple observations and consider a set $\set$ and
a group G acting on $\set$. For any such group action there exists some structure on $\set$ such that
G is a subgroup of the full automorphism group of this structure. If the set is finite then the group is
actually the full group of automorphisms of the space $\set$. The algebra of scalar difference operators on 
$\set$ will
be the skew group algebra $A=\fl[G]$ where $\fl$ will be the algebra of $\f$-valued functions defined on $\set$.
Difference equations on $\set$ and their solutions must be invariant objects under the action of the
group G. If they are not invariant their description and solutions will depend on the
arbitrariness in the specification of the underlying space. The Klein Erlanger program in geometry
has shown that the building blocks of the geometry of a set with a group action are the invariants
of the group. Geometrical objects and their relations are constructed from invariants. In this way the
geometry will not depend on the arbitrariness of the underlying space. What we propose in this paper
is in the spirit of the Erlanger program in geometry. 

We propose that the building blocks of the theory of 
difference equations on a finite space with some structure are the invariants of the group of automorphisms
acting on the space. The algebra of difference operators will be the skew group algebra,A, of G and all main
notions in the theory of difference equations will be defined in terms of invariants. A linear difference
equations $\e$ will be an A-module, symmetries of $\e$ will be A-endomorphisms of $\e$. All conserved quantities
and structures of the equations will be invariant elements in the tensoralgebra of the equation $\e$. A special
role will be played by the equations corresponding to indecomposable and simple A-modules. 

In this paper we introduce a 
Categorical point of view on equations and solutions. The equations are objects in a full subcategory of the
Category of A-modules. Solutions of an equation are descriptions of the equation in terms of other equations. Only
descriptions that are invariants are allowed and this leads to the idea of a solution of an equation in a Category
of A-modules as a morphism between the given equation and some other equation. So solutions are morphisms in the
Category. Solving an equation thus means to find the G-invariant descriptions or morphisms between the given
equation and all other equations in the Category. In this way symmetries are special types of solutions, they  are
descriptions of an equation in terms of itself. Simple equations play a special role in that they can only
be described in terms of themselves. They play the role of atoms in our category of equations. In the semisimple
situation all equations are sums of simple equations so the description of a given equation in terms of
simple equations in fact given a complete description of the equation. In a more general situation we also
need descriptions in terms of indecomposable equations in order to give a complete description of a given equation.
The indecomposable equations that are not simple are closely related to the notion of quantization.
The family of simple and indecomposable equations is determined by the group of symmetries of the underlying space
so this group determines the type of solutions that are needed to solve any equation in the Category.
Note that from this point of view a solution is a relational concept. It does not belong to one
object but to a pair of objects.

We will in this paper develope the theory for a class of equations we call finite type. These are analogs of the
finite type or Frobenius equations in the theory of differential equations. Note that if the set $\set$ is finite 
then all equations are of finite type. In an upcoming work the theory will be developed for a much wider class
of equations.

There exists currently several geometric-algebraic approaches to the study of difference and differential equations;
the differential algebra approach of Ritt \cite{ritt} and Kolchin \cite{kol} and the description through use of
jet bundles and D-modules \cite{lyc} just to mention two. Our approach does not belong directly to any of these directions.
It is however somewhat related to the approach in \cite{lyc} and the difference algebra approach in \cite{kon}.

\section{The main notions in the theory of finite type difference equations on a set}

Let $\set$ be a set and let G be a group acting on $\set$. We will assume that the
action of G is from the left and is faithful and transitive so that G is acting as a transitive group
of permutations on $\set$. It is well know from the theory of permutation groups that 
there exists a finite set of relations on $\set$ such that the group G is included in the full
group of symmetries of these relations. A space is a set with some
structure defined. Any group acting on a set can thus be thought of a the symmetry group
of a space. Examples of such spaces are graphs, lattices, finite projective
spaces, finite linear spaces etc.
\begin{example}
Let $\gr$ be the cyclic graph with vertex set $\{x_1,x_2,\cdots,x_n\}$ and edge set $\{\{x_1,x_2\},\{x_2,x_3\},
\cdots,\{x_n,x_1\}\}$. This graph can be considered to be a discrete approximation to
the circle $S^1$. The group of symmetries of this finite space is the dihedral group, $D_{2n}$. It has
two generators $t$ and $s$, where $t$ is reflection around $x_1$ and $s$ is left translation.
The symmetry group has $2n$ elements and presentation in terms of generators and
relations in the following form, $D_{2n}=\langle s,t\mid s^n=1,t^2=1,tst=s^{n-1}\rangle$.
\end{example}
We will now define the main notions in the theory of difference equations on a space S as
the invariants associated to the group of symmetries of S.

\subsection{The algebra of G-difference operators}

Let $\f$ be a field and let $\fl$ be the $\f$-algebra of $\f$ valued functions on S.
Let G be the symmetry group of the space S. The the left action of G on S can be lifted
to a left  action of G on $\fl$ in the natural way
\begin{equation*}
(g(f))(x)=f(s^{-1}x)\quad\forall g\in G
\end{equation*}
The skew group algebra of G over $\fl$ is the set of finite formal linear combinations
of elements of G with coefficients in $\fl$. Addition and multiplication by elements
in $\f$ is defined componentwise and multiplication is defined by
\begin{equation*}
(fg)(hg')=(fg(h))gg'
\end{equation*}
We now define the first basic notion in our theory of difference equations.

\begin{definition}
	$\flg$ is the algebra of $\emph{G-difference operators}$ on the set S.
\end{definition}

\begin{notation}
We will from now on use the notation $k=\fl$ and $A=\flg$.
\end{notation}

\begin{example}
Let $\set$ be the cyclic graph $\gr$ with n elements. We have seen that the symmetry group of
$\gr$ is the dihedral group $D_{2n}$ with generators beeing left translation s  and reflection t.
The algebra of $D_{2n}$-difference operators consists of formal linear combinations of $\f$-valued
functions and elements of $D_{2n}$. The algebra A contains the usual difference operators from
the calculus of differences whose continum limit corresponds
to the usual ordinary differential operators. But it also contains operators involving the reflection t.
These operators will in the continum limit correspond to differential-delay equations.
\end{example}

\subsection{Linear G-difference equations of finite type and solutions}

Let A be the algebra of G-difference operators on a space $\set$. Let $\e$
be a finitely generated module over k. If not otherwise noted finitely generated
means finitely generated over k. Assume that G acts on $\e$ on the left, $g(fe)=g(f)ge$. Then
$\e$ is a left A-module with the natural action of the skew group algebra A on $\e$.
In this way $\e$ can be considered to be an invariant for the symmetry group G of the
underlying space. We will consider only left A-modules that can be given a geometrical 
interpretation. Define $\mu_x\sub k$ by

\begin{equation*}
\mu_x=\{f:\set\rightarrow\f\mid\;f(x)=0\}
\end{equation*}

The subsets $\mu_x$ are clearly ideals in k. They are in fact maximal ideals 

\begin{proposition}
$\mu_x$ is a maximal ideal in $\fl$.

\end{proposition}
\begin{proof}
	$\mu_x$ is clearly an ideal in $\fl$. Let $J$ be an ideal in $\fl$ and assume that
	$\mu_x\sub J\sub\fl$ $\Rightarrow$ $\exists$ $j\in J$ such that $j\notin\mu_x$ 
$\Rightarrow$	 $j(x)\ne 0$ $\Rightarrow$ $\frac{j}{j(x)}\in J$. But then $(\frac{j}{j(x)}-1)(x)=
\frac{j(x)}{j(x)}-1=0$	 $\Rightarrow$ $\frac{j}{j(x)}-1\in\mu_x$. But $\mu_x\sub J$ 
$\Rightarrow$ $\frac{j}{j(x)}-1\in J$	 $\Rightarrow$ $\frac{j}{j(x)}-(\frac{j}{j(x)}-1)=1\in J$ 
$\Rightarrow$ $J=\fl$. This is a contradiction	 so $J=\mu_x$ and $\mu_x$ is maximal.
\end{proof}

For each $x\in \set$ we have a submodule $\mu_x\e$ since $\mu_x$ is an ideal. We will only
consider left A-modules that have no invisible elements \cite{lyc}.

\begin{definition}
$\e$ is a \emph{geometric} left A-module iff

\begin{equation}
\bigcap_{x\in\set}\mu_x\e=0
\end{equation}
\end{definition}

$\e$ beeing a left A-module means that we have an action of the algebra of G-difference operators
on $\e$. We are now ready to define the second main notion in our theory.

\begin{definition}
A linear G-difference equation of finite type is a geometric left A-module
that is finitely generated over k.
\end{definition}

We will use the term GF-difference equations for the equations defined in the previous
definition.
In general the structure of a GF-difference equation is investigated by comparing it to other
equations. An equation will be considered to be understood
only if its relations to all other equations are known.  This is the Categorical point of view.
Relations between equations are A-morphisms so an equation $\e$ is understood or $\emph{solved}$ if 
$Hom_A(\e,\fc)$ is known for all GF-difference equations $\fc$. Let us formalize this in a definition
\begin{definition}
Let $\e$ be any GF-difference equation. Then a solution of $\e$ of type $\fc$, where $\fc$ is a GF-difference
equation,  is a A-module morphism $\phi\in Hom_A(\e,\fc)$
\end{definition}
Using this definition we can now say that a GF-difference equation is solved if we know all solutions
of the equation. We will introduce two special types of solutions that will play a central
role in our theory. A GF-difference equation is indecomposable if it can not be written as a direct
sum of two GF-difference equations. Our first special type of solution is the following

\begin{definition}
Let $\e$ be any GF-difference equation and let $\s$ be a indecomposable equation. Then a indecomposable 
solution of $\e$ of type $\s$ is an element of $Hom_A(\e,\s)$.
\end{definition}

The second special type of solution is symmetries. These are relations that
describe the equation in terms of itself so we define.

\begin{definition}
Let $\e$ be any GF-difference equation. Then a symmetry of $\e$ is a A-morphism of $\e$
to itself.
\end{definition}

So a symmetry of $\e$ is an element of $End_A(\e)$. If $f\in Hom_A(\e,\fc)$
is any solution $\e$ of type $\fc$ and $\phi$ is a symmetry of $\e$ then $\phi^*(f)=\phi\circ f$
is also a solution of type $\fc$. So symmetries map solutions of some type to solutions
of the same type.  The problem of solving an equation is closely linked to the module structure of
the equation and we will now start to develope the structure theory for GF-difference equations.

\section{The structure of the Category of GF-difference equations}

Let $\GFE$ be a category \cite{lane} whose objects are GF-difference equations and morphisms are
A-module morphisms.
\begin{definition}
$\GFE$ is the category of GF-difference equations
\end{definition}
A complete description of the structure of the category $\GFE$ is the same as knowing all
solutions to all GF-difference equations. This is in general an enormously complicated
problem. We will in this section describe what can be said in general about the structure
of the category $\GFE$ without placing any restrictions on the set $\set$ or the group G.
We will start our investigation of the structure of $\GFE$ by investigating the closure of
the set of GF-difference with respect to the usual linear
algebra operations like direct sum, tensor product etc. These
operations preserve the set of modules that are finitely generated over k. They
also preserve the property of beeing geometric as we will now see.

\subsection{The algebra of GF-difference equations}

Let $\e_1,\e_2$ be two GF-difference equations. Then $\e_1\ds\e_2$ is a 
finitely generated left A-module with the operations

\begin{align*}
f(e_1,e_2)&=(fe_1,fe_2)\quad\text{$\forall f\in k$}\\
g(e_1,e_2)&=(ge_1,ge_2)\quad\text{$\forall g\in G$}.
\end{align*}

\begin{proposition}
The direct sum of GF-difference equations $\e_1\ds\e_2$ is a GF-difference equation.
\end{proposition}
\begin{proof}
We know that $\bigcap_{x\in\set}\mu_x\e_1=0,\;\bigcap_{x\in\set}\mu_x\e_2=0$ and we have by definition
that $\mu_x\e=\{\sum_{i}f_ie_i\mid\;f_i\in\mu_x,e_i\in\e\}$ $\imp$ $\mu_x(\e_1\ds\e_2)=
\{\sum_{i}f_i(e_i^1+e_i^2)\mid\;f_i\in\mu_x,e_i^1\in\e_1,e_i^2\in\e_2\}$. But $\sum_if_i(e_i^1
+e_i^2)=\sum_if_ie_i^1+\sum_if_ie_i^2\in\mu_x\e_1\ds\mu_x\e_2$. So $\mu_x(\e_1\ds
\e_2)\sub\mu_x\e_1\ds\mu_x\e_2$. $\imp$ $\bigcap_{x\in\set}\mu_x(\e_1\ds\e_2)
\sub\bigcap_{x\in\set}(\mu_x\e_1\ds\mu_x\e_2)\sub\bigcap_{x\in\set}\mu_x\e_1\ds
\bigcap_{x\in\set} \mu_x\e_2=0+0=0$. So the direct sum is a finitely generated geometric left A-module.
\end{proof}

Let $\e_1,\e_2$ be GF-difference equations. From this it follows that they are left k-modules since
$k\sub A$ as algebras. The algebra k is abelian so $\e_1\tpk\e_2$ is a well
defined finitely generated k-module. Define a $G$ action on the tensorproduct module by $g(e_1\tpk e_2)=ge_1
\tpk ge_2$. With this action we have

\begin{proposition}
$\e_1\tpk\e_2$ is a A-module with the given $G$ action.
\end{proposition}
\begin{proof}
\begin{align*}
g(fe_1\tpk e_2)&=g(fe_1)\tpk ge_2=g(f)(ge_1)\tpk ge_2\\
&=ge_1\tpk g(f)(ge_2)=ge_1\tpk g(fe_2)\\
&=g(e_1\tpk fe_2).\\
(g_1g_2)(e_1\tpk e_2)&=(g_1g_2)e_1\tpk (g_1g_2)e_2=g_1(g_2e_1)\tpk 
g_1(g_2e_2)\\
&=g_1(g_2e_1\tpk g_2e_2)=g_1(g_2(e_1\tpk e_2)).\\
(gf)(e_1\tpk e_2)&=g(fe_1\tpk e_2)=g(fe_1)\tpk ge_2\\
&=(g(f)g)e_1\tpk ge_2=(g(f))(ge_1\tpk ge_2)=(g(f)g)(e_1\tpk e_2).
\end{align*}
\end{proof}

\begin{proposition}
Assume that $\e_1,\e_2$ are two GF-difference equations. Then $\e_1\tpk\e_2$ is also
a GF-difference equation.
\end{proposition}
\begin{proof}
We observe that $\mu_x(\e_1\tpk\e_2)\sub\mu_x\e_1\tpk\e_2$. So we have $\sn_{x\in\set}\mu_x(\e
_1\tpk\e_2)\sub\sn_{x\in\set}(\mu_x\e_1\tpk\e_2)\sub(\sn_{x\in\set}\mu_x\e_1)\tpk\e_2=0$.
\end{proof}

Let $\e_1,\e_2$ be two GF-difference equations. Then $Hom_k(\e_1,\e_2)$ is a finitely generated
left k-module with the natural action of k

\begin{equation*}
(f\phi)(e)=f(\phi(e))\quad\text{$\forall f\in k$}
\end{equation*}

Define an action of $G$ on $Hom_k(\e_1,\e_2)$ by

\begin{equation*}
(g\phi)(e)=g(\phi(g^{-1}e))\quad\text{$\forall g\in G$}
\end{equation*}

\begin{proposition}
$Hom_k(\e_1,\e_2)$ is a left A-module with the given action of k and $G$.
\end{proposition}
\begin{proof}
\begin{align*}
(g\phi)(fe)&=g(\phi(g^{-1}fe))=g(\phi((g^{-1})fg^{-1}e))=g((g^{-1})f\phi(g^{-1}e))\\
&=g(g^{-1})f(g\phi(g^{-1}e))=f((g\phi)(e))\\
((g_1g_2)\phi)(e)&=(g_1g_2)\phi((g_1g_2)^{-1}e)=g_1(g_2\phi(g_2^{-1}(g_1^{-1}e)))\\
&=g_1((g_2\phi)(g_1^{-1}e))=(g_1(g_2\phi))(e)\\
((gf)(\phi))(e)&=(g(f\phi))(e)=g((f\phi)(g^{-1}e))=g(f(\phi(g^{-1}e)))\\
&=g(f)(g(\phi(g^{-1}e)))=g(f)((g\phi)(e))=(g(f)(g\phi))(e)=((g(f)g)(\phi))(e).
\end{align*}
\end{proof}

\begin{proposition}
$Hom_k(\e_1,\e_2)$ is GF-difference equation.
\end{proposition}
\begin{proof}
Let $\psi\in\mu_xHom_k(\e_1,\e_2)$ $\imp$ $\psi=\sum_if_i\phi_i$ with $f_i\in\mu_x,\phi_i\in
 Hom_k(\e_1,\e_2)$ $\imp$ $\psi(e)=\sum_if_i(\phi(e))\in\mu_x\e_2$ so $\psi(e)\in\mu_x\e_2$
for all $e\in\e_1$ $\imp$ $\psi\in Hom_k(\e_1,\mu_x\e_2)$. So we have $\mu_x Hom_k(\e_1,\e_2)\sub
 Hom_k(\e_1,\mu_x\e_2)$. But then $\sn_{x\in\set}\mu_x Hom(\e_1,\e_2)\sub \sn_{x\in\set}
Hom_k(\e_1,\mu_x\e_2)\sub Hom_k(\e_1,\sn_{x\in\set}\mu_x\e_2)=0$
\end{proof}

As a special case of the last proposition we have

\begin{corollary}
Let $\e$ be a GF-difference equation. Then the dual $\e^*$ is also a GF-difference equation.
\end{corollary}

Let us next consider the case of quotients. Assume that $\e$ is a GF-difference equation and let
$\e'\sub\e$ be a submodule of $\e$.

\begin{proposition}
$\e/\e'$ is GF-difference equation.
\end{proposition}
\begin{proof}
Since $\e/\e'$ is finitely generated we only need to prove that it is geometric.
Let $h\in\mu_x(\e/\e')$ then  $h=\sum_if_i[e_i]$ where $f_i\in\mu_x$ and $[e_i]\in\e/\e'$ are 
the equivalent classes of elements in $\e$. So $h=\sum_if_i[e_i]=[\sum_if_ie_i]$ and we can 
conclude $\mu_x(\e/\e')\sub\mu_x\e/\e'$. But then we have $\sn_{x\in\set}\mu_x(\e/\e')\sub 
\sn_{x\in\set}\mu_x\e/\e'\sub(\sn_{x\in\set}\mu_x\e)/\e'=0$.
\end{proof}

We already know that tensor products and direct sums of GF-difference equations are
GF-difference equations. This implies that the tensor algebra $T\e$ of a GF-difference equation
is a GF-difference equation. The modules $S^n\e$ and $\wedge^n\e$ are 
factors of the tensor algebra of $\e$ so we have the following result.

\begin{corollary}
Let $\e$ be a GF-difference equation. Then $S^n\e$ and $\wedge^n\e$ are GF-difference equations.
\end{corollary}

\subsection{GF-difference equations as modules of sections in vectorbundles}

We have seen that the category of GF-difference equations is closed with respect to quotients,
$\ds$,$\tpk$,$Hom_k$,
$\wedge_k$ and $S^n_k$. These modules can be given an interpretation as modules of sections
in vectorbundles over the set $\set$. 

 Let $\e$ be a GF-difference equation. Then in 
particular $\e$ is a k-module and $\mu_x\e\sub\e$ is a k submodule of $\e$. Let 
$E_x=\e/\mu_x\e$.
Then $E_x$ is A-module over $k/\mu_x\approx\f$ and therefore is a $\f$ vectorspace. Denote the 
elements of $E_x$ by $[e]_x$ where $[e]_x=[e']_x$ only if $e-e'\in\mu_x\e$. Let the bundle $B$ over
$\set$ be defined by

\begin{equation*}
B=\un_{x\in\set}(x,E_x)
\end{equation*}

where the projection,$\pi:B\ft\set$, in the bundle is projection on the first component. Let $\gm(B)$ be the set of
sections in the bundle $B$. This set is a module over k through pointwise addition and 
multiplication by functions in k.

For each element in $G$ define a bundle map in the bundle $B$ through

\begin{equation}
g(x,[e]_x)=(gx,[ge]_{gx}).\label{act1}
\end{equation}

This set of bundle maps in fact defines an action of $G$ on the bundle $B$.

\begin{proposition}
Bbundle map \ref{act1} is well defined and determines an action of $G$ on the 
bundle $B$.
\end{proposition}
\begin{proof}
Assume $[e]_x=[e']_x$. Then $e-e'\in\mu_x\e$ and so $ge-ge'\in g(\mu_x\e)$. Let $\tilde{e}\in\mu_x\e$,
then $\tilde{e}=\sum_if_ie_i$ where $e_i\in\e$ and $f_i\in\mu_x$ and so $g\tilde{e}=\sum_ig(f)_i
ge_i$. But $g(f)_i(gx)=f_i(x)=0$ so $g(f_i)\in\mu_{gx}$. Then it follows that $ge-ge'\in\mu_{gx}\e$ and 
we can conclude that $[ge]_{gx}=[ge']_{gx}$ so the map is well defined. Using the definition of the 
bundle map we have
\begin{align*}
(g_1g_2)(x,[e]_x)&=((g_1g_2)x,[(g_1g_2)e]_{(g_1g_2)x})=(g_1(g_2x),[g_1(g_2e)]_{g_1(g_2x)})\\
&=g_1(g_2x,[g_2e]_{g_2x})=g_1(g_2(x,[e]_x))
\end{align*}

so the bundle map defines an action of $G$ on the bundle $B$.
\end{proof}

\begin{corollary}
$B$ is a vectorbundle,that is $dimE_x$ is constant.
\end{corollary}
\begin{proof}
Let $x\in\set$ be a fixed point in the set $\set$. The group acts transitively on the set  $\set$
so for any $y\in\set$ there exists a $g\in G$ such that $gx=y$.This element induces a map 
$\phi_g:E_x\ft E_y$ defined by $\phi_g([e]_x)=[ge]_{gx}$. This map is linear and has an inverse
$\phi_{g^{-1}}$. We can therefore conclude that all fibers $E_y$ of the bundle $B$ have the same 
dimension so $B$ is a vectorbundle.
\end{proof}

We now induce an action of $G$ on $\gm(B)$ defining
\begin{equation}
(gs)(x)=g(s(g^{-1}x))\label{act2}
\end{equation}

\begin{proposition}
Action \ref{act2} gives $\gm(B)$ the structure of an A-module.
\end{proposition}
\begin{proof}
\begin{align*}
((gf)(s))(x)&=(g(fs))(x)=g((fs)(g^{-1}x))=g(g(f)(x)s(g^{-1}x))\\
&=g(f)(x)(gs)(x)=(g(f)g)(s)(x)\\
((g_1g_2)s)(x)&=(g_1g_2)s((g_1g_2)^{-1}x)=g_1(g_2s(g_2^{-1}(g_1^{-1}x)))\\
&=g_1((g_2s)(g_1^{-1}))=(g_1(g_2s))(x)
\end{align*}
\end{proof}

\begin{proposition}
$\e\iso\gm(B)$ as A-modules.
\end{proposition}
\begin{proof}
Let $e\in\e$. Define $\phi(e)\in\gm(B)$ by $\phi(e)(x)=[e]_x$. We clearly have 
$\phi:\e\ft\gm(B)$.
\begin{align*}
\phi(e+f')(x)&=[e+f']_x=[e]_x+[f']_x=\phi(e)(x)+\phi(f')(x)=(\phi(e)+\phi(f'))(x),\\
\phi(fe)(x)&=[fe]_x=f(x)[e]_x=f(x)\phi(e)(x)=(f\phi(e))(x).
\end{align*}
We have used that fact that $z[e]_x=[fe]_x$ where $f$ is any function such that $f(x)=z$.
This is well defined because if $[e]_x=[e']_x$ and $f(x)=z,f'(x)=z$ then $e'=e+h,f'=f+g$
where $h\in\mu_x\e$ and $g\in\mu_x$. But then $fe-f'e'=fe-(f+g)(e+h)=-ge-fh-gh\in\mu_x\e$.
So $[fe]_x=[f'e']_x$.

We have now proved that $\phi$ is a k-module morphism. It is also a A-module morphism

\begin{align*}
\phi(ge)(x)&=[ge]_x=[ge]_{g(g^{-1}x)}=g([e]_{g^{-1}x})=g(\phi(e)(g^{-1}x))=(g\phi(e))(x).
\end{align*}

Assume $\phi(e)=\phi(e')$. Then $[e]_x=[e']_x$ so $[e-e']_x=0\quad\forall
 x\in\set$. But this means that $e-e'\in\sn_{x\in\set}\mu_x\e$. We can therefore conclude that
$e=e'$ because $\e$ is geometric. So $\phi$ is injective.

 For each $y\in\set$ let $\Pi_y:\e\ft E_y$ be
the canonical projection. $B$ is a vectorbundle so that the dimension n of each fibre as a $\f$-vector space
is the same. Let $\{e_i\}_{i=1}^m$ be a set of generators for $\e$. Then $\Pi_y(\{e_i\}_{i=1}^m)$ generates
$E_y$ for all $y\in\set$ so at each point at least one subset of say n elements of $\{e_i\}_{i=1}^n$ form a
basis for $E_y$ after projection by $\Pi_y$. There are only finitely many subsets of n elements from the
set of m generators. Enumerate these subsets
\begin{align*}
B^i=\{e_{l(i,k)}\}_{k=1}^n \quad \text{$i=1\cdots r$}.
\end{align*}

Here $l(i,k)$ is an index function. Put $\set_1=\set$ and define subsets $V_i\sub\set$ recursively
\begin{align*}
V_i&=\{y\in\set_i\mid \Pi_y(B^i) \quad\text{is a basis of $E_y$}\},\\
\set_{i+1}&=\set_i-V_i.
\end{align*}
This gives us a finite set of nonempty subsets $\{V_i\}_{i=1}^p$ such that $V_i\sn V_j=\emptyset$ for $i\ne j$,
$\set=\un_{i=1}^pV_i$ and $\Pi_yB^i=\{[e_{l(i,k)}]_y\}_{k=1}^n$ is a basis for $E_y$ for all $y\in\set$. Let
$\delta_{V_i}$ be the characteristic function for $V_i$. Then $\delta_{V_i}\in k$ and $\sum_i\delta_{V_i}=1$
Let $T_i=\delta_{V_i}\gm(B)$. Then $T_i$ is a k-submodule of $\gm(B)$ and $T_i$ has k-basis 
$\{\delta_{V_i}\phi(e_{l(i,k)})\}_{k=1}^n$. Let $s\in\gm(B)$ be any section. Then we have
\begin{align*}
s&=(\sum_i\delta_{V_i})(s)=\sum_i\delta_{V_i}s=\sum_i\sum_kf_{ik}\delta_{V_i}\phi(e_{l(i,k)})\\
&=\phi(\sum_i\sum_kf_{ik}\delta_{V_i}e_{l(i,k)}).
\end{align*}
But $e=\sum_i\sum_kf_{ik}\delta_{V_i}e_{l(i,k)}\in\e$ so $\phi$ is surjective.

\end{proof}

This result show that the category $\GFE$ is equivalent to the category of modules of sections in
vectorbundles $\gm(B)$ over $\set$ where the action of G is defined through $\ref{act1}$ and $\ref{act2}$.

\begin{proposition}
Let $\e$ be a GF-difference equation. Then $\e$ is free and finitedimensional as a k-module.
\end{proposition}
\begin{proof}
Let $\{e_i^x\}_{i=1}^n$ be a basis for $E_x$ over $\f$. The number $n$ exists and
is independent of x since we  all our bundles are finite dimensional and vector bundles
so that the dimension of all fibers are the same. Define sections $\{s_i\}_{i=1}^n$ by
$s_i(x)=e_i^x$. Assume that $\sum_{i=1}^nf_is_i=0$. Then $\sum_{i=1}^nf_i(x)e_i^x=0$ so
$f_i(x)=0$ for all $x\in\set$ and all $i$. This implies that $f_i=0$ for all $i$ and we conclude
that $\{s_i\}_{i=1}^n$ is a linearly independent set over k. Let $s\in\gm(B)$, then 
$s(x)=e_x\in E_x$ so there exists complex numbers $\{c_i^x\}_{i=1}^n$ such that
 $s(x)=\sum_{i=1}^nc_i^xe_i^x$. Define functions in k by $f_i(x)=c_i^x$, then 
$s=\sum_{i=1}^xf_is_i$ and $\{s_i\}_{i=1}^n$ is a spanning set.
\end{proof}

 We can use this result to
prove a standard isomorphism. Define a map $\phi:\e^*\times\fc\ft Hom_k(\e,\fc)$
by 
\begin{align*}
\phi(e^*,f)(e)=e^*(e)f
\end{align*}

\begin{proposition}
$\phi$ is k-bilinear.
\end{proposition}
\begin{proof}
\begin{align*}
\phi(e_1^*+e_2^*,f)(e)&=\phi((e_1+e_2)^*,f)(e)=(e_1+e_2)^*(e)f=e_1^*(e)f+e_2^*(e)f\\
&=\phi(e_1^*,f)(e)+\phi(e_2^*,f)(e)=(\phi(e_1^*,f)+\phi(e_2^*,f))(e)\\
\phi(e^*,f_1+f_2)(e)&=e^*(e)(f_1+f_2)=e^*(e)f_1+e^*(e)f_2\\
&=(\phi(e_1^*,f_1)+\phi(e^*,f_2))(e)\\
\phi(re^*,f)(e)&=(re^*)(e)f=(r(e^*(e)))f=r(e^*(e)f)\\
&=r(\phi(e^*,f)(e))=(r\phi(e^*,f))(e)
\end{align*}
\end{proof}

This proposition show that we have a well defined map $\phi:\e^*\tpk\fc\ft Hom_k(\e,\fc)$
defined by $\phi(e^*\tpk f)(e)=e^*(e)f$. 
\begin{proposition}
$\phi$ is an A-isomorphism
\end{proposition}
\begin{proof}
Let $\{e_i\}_{i=1}^n,\{f_i\}_{i=1}^m$ be basis over k for $\e$ and $\fc$. Let $\{e_i^*\}_{i=1}^n$ be the 
dual basis for $\e^*$. Then $\{e_i^*\tpk f_j\}$ is a basis for $\e^*\tpk\fc$ because the modules
are free over k. Let $v=\sum_{ij}a_{ij}e_i^*\tpk f_j$ and assume that $\phi(v)=0$. Then $\phi(v)(e_s)=0$
for all s and we have $\sum_ja_{sj}f_j=0$ so that $a_{ij}=0$ for all i and j because $\{f_j\}$ is 
a basis for $\fc$. So $\phi$ is injective. Let $F\in Hom_k(\e,\fc)$. Define the matrix $(F_{ij})$ by
$F(e_i)=\sum_jF_{ij}f_j$ and let $v=\sum_{ij}F_{ij}e_i^*\tpk f_j$. Then $\phi(v)=F$ so that $\phi$
is surjective. Finally we have
\begin{align*}
\phi(g(e^*\tpk f))(e)&=\phi(ge^*\tpk gf)(e)=(ge^*)(e)gf=(g(e^*(g^{-1}e)))gf\\
&=g(e^*(g^{-1}e)f)=g(\phi(e^*\tpk f)(g^{-1}e))=(g\phi(e^*\tpk f))(e)
\end{align*}
so that $\phi$ is an A-morphism.
\end{proof}

\begin{corollary}
$Hom_k(\e,\fc)\iso\e^*\tpk\fc$
\end{corollary}

We know that the category of GF-difference equations is closed with respect to the usual linear
algebra operations. Since we have proved that any GF-difference equation is isomorphic to A-module of 
sections in
a vectorbundle over $\set$ it is evident that all such linear algebra operations must 
reduce to operations on the corresponding vectorbundles. The following series of
propositions show that the correspondence is as nice as one would expect.

Let $\e_1,\e_2$ be two GF-difference equations.. Then $\e_1\iso\gm(B_1)$ and $\e_2\iso\gm(B_2)$
where $B_1$, $B_2$ are vectorbundles

\begin{align*}
B_1=\un_{x\in\set}(x,E_x^1)\\
B_2=\un_{x\in\set}(x,E_x^2)
\end{align*}

We then have

\begin{proposition}
$\e_1\ds\e_2\iso\gm(B_1\ds B_2)$.
\end{proposition}
\begin{proof}
Define a map $\phi:\e_2\ds\e_2\ft\gm(B_1\ds B_2)$ by
\begin{equation*}
\phi(s_1,s_2)(x)=(s_1(x),s_2(x))
\end{equation*}
where we identify the GF-difference equations  $\e_1,\e_2$ with their corresponding
modules of sections.

Assume $\phi(s_1,s_2)=\phi(s_1',s_2')$. Then $(s_1(x),s_2(x))=(s_1'(x),s_2'(x))$
and so  $s_1(x)=s_1'(x),s_2(x)=s_2'(x)$ for all $x\in\set$. But this implies that
$(s_1,s_2)=(s_1',s_2')$ and $\phi$ is injective.

Let $s\in\e_1\ds\e_2$. Then $s(x)\in E_x^1\ds E_x^2$ for all $x\in\set$. Define 
$s_1(x)=\pi_1\circ s(x)$ and $s_2(x)=\pi_2\circ s(x)$ where $\pi_1:E_x^1\ds E_x^2\ft E_x^1$
and $\pi_2:E_x^1\ds E_x^2\ft E_x^2$ are the projections on the first and second 
factor. But then $(s_1,s_2)\in\e_1\ds\e_2$ and evidently 
\begin{equation*}
\phi(s_1,s_2)(x)=(s_1(x),s_2(x))=s(x)
\end{equation*}
so $\phi$ is surjective. Furthermore we have
\begin{align*}
\phi(f(s_1,s_2))(x)&=\phi(fs_1,fs_2)(x)=((fs_1)(x),(fs_2)(x))\\
&=(f(x)s_1(x),f(x)s_2(x))=f(x)(s_1(x),s_2(x))=(f\phi(s_1,s_2))(x),\\
\phi(g(s_1,s_2))(x)&=\phi(gs_1,gs_2)(x)=((gs_1)(x),(gs_2)(x))\\
&=(g(s_1(g^{-1}x)),g(s_2(g^{-1}x)))=g((s_1,s_2)(g^{-1}x))=(g\phi(s_1,s_2))(x).
\end{align*}

So $\phi$ is a left A-module morphism and the proof is complete.
\end{proof}

We have seen that the k-tensor product of GF-difference equations is a GF-difference
equation with the action of k and $G$ defined by
\begin{align*}
f(s_1\tpk s_2)&=fs_1\tpk s_2\quad\text{$\forall f\in k$},\\
g(s_1\tpk s_2)&=gs_1\tpk gs_2\quad\text{$\forall g\in G$}.
\end{align*}

Using the vectorbundles $B_1$ and $B_2$ corresponding to $\e_1$ and $\e_2$ we define a
new vectorbundle $B_1\tpc B_2$ by
\begin{equation*}
B_1\tpc B_2=\un_{x\in\set}(x,E_x^1\tpc E_x^2).
\end{equation*}

Let $\gm(B_1\tpc B_2)$ be the set of sections in the vectorbundle $B_1\tpc B_2$. This set
is a k-module through pointwise addition and multiplication by elements of k. It is also
a left A-module through the action
\begin{equation*}
(gs)(x)=g(s(g^{-1}x))
\end{equation*}
where
\begin{equation*}
g(x,[e_1]\tpc [e_2])=(gx,[ge_1]_{gx}\tpc [ge_2]_{gx}).
\end{equation*}
We then have the following result

\begin{proposition}
$\e_1\tpk\e_2\iso\gm(B_1\tpc B_2)$.
\end{proposition}
\begin{proof}
Define a map $\tilde{\phi}:\e_1\times\e_2\ft\gm(B_1\tpc B_2)$ by
\begin{equation*}
\tilde{\phi}(s_1,s_2)(x)=s_1(x)\tpc s_2(x).
\end{equation*}
We have
\begin{align*}
\tilde{\phi}(s_1+s_1,s_2)&=(s_1+s_1)(x)\tpc s_2(x)=(s_1(x)+s_1(x))\tpc s_2(x)\\
&=s_1(x)\tpc s_2(x)+s_1(x)\tpc s_2(x)=(\phi(s_1,s_2)+\phi(s_1,s_2))(x),\\
\tilde{\phi}(fs_1,s_2)(x)&=(fs_1)(x)\tpc s_2(x)=(f(x)s_1(x))\tpc s_2(x)\\
&=s_1(x)\tpc (f(x)s_2(x))=s_1(x)\tpc (fs_2)(x)=\phi(s_1,fs_2)(x).
\end{align*}
So $\tilde{\phi}$ is k-bilinear and therefore induces a unique map 
$\phi:\e_1\tpk\e_2\ft\gm(B_1\tpc B_2)$ where 
\begin{equation*}
\phi(s_1\tpk s_2)(x)=s_1(x)\tpc s_2(x).
\end{equation*}
Let $\{s_i^1\}_{i=1}^{n}$ and $\{s_i^2\}_{i=1}^{n}$ be bases for $\e_1$ and $\e_2$ as k-modules.
These bases exists because the modules are free as modules over k. Let $s\in\e_1\tpk \e_2$ then
\begin{equation*}
s=\sum_{ij}f_{ij}s_i^1\tpk s_j^2.
\end{equation*}
Assume that $\phi(s)=0$. This implies that
\begin{equation*}
s=\sum_{ij}f_{ij}(x)s_i^1(x)\tpc s_j^2(x)=0\quad\text{$\forall x\in\set$}.
\end{equation*}
But then $f_{ij}(x)=0\quad\forall x\in\set$ and so $f_{ij}=0$ and as a consequence $s=0$. So
$\phi$ is injective.
Let $s\in\gm(B_1\tpc B_2)$. This implies that $s(x)\in B_1\tpc B_2$ so there exists elements of $\f$
$z_{ij}^x$ such that 
\begin{equation*}
s(x)=\sum_{ij}z_{ij}^xs_1(x)\tpc s_2(x).
\end{equation*}
Define elements $f_{ij}\in$  k by $f_{ij}(x)=z_{ij}^x$ and define $h\in\e_1\tpk\e_2$ by
\begin{equation*}
h=\sum_{ij}f_{ij}s_i^1\tpk s_j^2.
\end{equation*}
Then we evidently have $\phi(h)=s$ and $\phi$ is surjective.

We already know that $\phi$ is k-linear by definition of tensor product. But we also have
\begin{align*}
\phi(g(s_1\tpk s_2))(x)&=\phi(gs_1\tpk gs_2)(x)=(gs_1)(x)\tpc (gs_2)(x)\\
&=g(s_1(g^{-1}x))\tpc g(s_2(g^{-1}x))=g((s_1\tpk s_2)(g^{-1}x))\\
&=g(\phi(s_1,s_2)(g^{-1}x))=(g\phi(s_1\tpk s_2))(x)
\end{align*}

so $\phi$ is a A-module morphism and the proof is complete
\end{proof}

Now let $\e\iso\gm(B)$ be a GF-difference equation with corresponding vectorbundle $B$. Let
$\e'\sub\e$ be a subequation. Define
\begin{equation*}
V_x=\{[e]_x\mid e\in\e'\}.
\end{equation*}
Then $V_x\sub E_x$ is a subspace of $E_x$ for each $x\in\set$ and the dimension is independent
of x. Define a vectorbundle $B'$ by
\begin{equation*}
B'=\un_{x\in\set}(x,V_x).
\end{equation*}
Then $B'$ is evidently a subvectorbundle of $B$ and we have by construction that $\e'\iso\gm(B')$.
Let $E_x/V_x$ be the factor space. Its dimension is independent of x and we can form the 
vectorbundle
\begin{equation*}
B/B'=\un_{x\in\set}(x,E_x/V_x).
\end{equation*}

Denote the elements of $E_x/V_x$ by $[v_x]_{V_x}$. We define an action by elements in $G$ by
\begin{equation*}
g([v_x]_{V_x})=[g(v_x)]_{V_{gx}}.
\end{equation*}

This action is well defined and we use it to induce an action of $G$ on $\gm(B/B')$ in the usual way.

\bpr
$\e/\e'\iso\gm(B/B')$
\epr
\bp
Define a map $\phi:\gm(B)/\gm(B')\ft\gm(B/B')$ by $\phi([s])(x)=[s(x)]_{V_x}$. Then $\phi$ is 
well defined  because if $[s]=[s']$ then $s-s'\in\gm(B')$ and therefore $s(x)-s'(x)\in V_x$. So
$[s(x)]_{V_x}=[s'(x)]_{V_x}$ and therefore $\phi([s])=\phi([s'])$.
Furthermore we have
\begin{align*}
\phi([s]+[s'])(x)&=\phi([s+s'])(x)=[(s+s')(x)]_{V_x}\\
&=[s(x)+s'(x)]_{V_x}=[s(x)]_{V_x}+[s'(x)]_{V_x}=(\phi([s])+\phi([s']))(x),\\
\phi(f[s])(x)&=\phi([fs])(x)=[(fs)(x)]_{V_x}=[f(x)s(x)]_{V_x}\\
&=f(x)[s(x)]_{V_x}=f(x)\phi([s])(x)=(f\phi([s]))(x),\\
\phi(g[s])(x)&=\phi([gs])(x)=[(gs)(x)]_{V_x}=[g(s(g^{-1}x))]_{V_x}\\
&=[g(s(g^{-1}x))]_{V_{g(g^{-1}x)}}=g([s(g^{-1}x)]_{V_{g^{-1}x}})\\
&=g(\phi([s])(g^{-1}x)=(g\phi([s]))(x).
\end{align*}
So we can conclude that $\phi$ is a A-module morphism.
Assume that $\phi([s])=\phi([s'])$. Then $\phi([s])(x)=\phi([s'])(x)$ for all $x\in\set$.
But this is the same as $[s(x)]_{V_x}=[s'(x)]_{V_x}$ so $s(x)-s'(x)\in V_x$. This implies that
$s-s'\in\gm(B')$ so by definition $[s]=[s']$ and $\phi$ is injective.
Let $\gamma\in\gm(B/B')$, then $\gamma(x)=[v_x]_{V_x}$. Define 
$s\in\gm(B)$ by $s(x)=v_x$. Then clearly $\phi([s])(y)=[s(y)]_{V_y}$ so $\phi([s])=\gamma$ and
$\phi$ is surjective.
\ep
Since $\wedge^n\e$ and $S^n\e$ are factor bundles of the tensor algebra $T\e$, it follows from the 
previous proposition that

\begin{corollary}
\begin{align*}
\wedge^n\gm(B)=\gm(\wedge^n B),\\
S^n\gm(B)=\gm(S^nB).
\end{align*}
\end{corollary}

Let $\e_1$ and $\e_2$ be GF-difference equations. We have proved that $Hom_k(\e_1,\e_2)$ is
a GF-difference equation with the actions
\begin{align*}
(f\phi)(s_1)&=f(\phi(s_1))\quad\text{$\forall f\in k$},\\
(g\phi)(s_1)&=g(\phi(g^{-1}s_1)).
\end{align*}
We know that $\e_1\iso\gm(B_1)$, $\e_2\iso\gm(B_2)$ where $B_1=\un_{x\in\set}(x,E_x^1)$ and
$B_2=\un_{x\in\set}(x,E_x^2)$ are vectorbundles. Let $Hom_{\f}(E_x^1,E_x^2)$ be
the set of $\f$-linear maps from $E_x^1$ to $E_x^2$. These have all the same dimension 
and we can form the vector bundle
\begin{equation*}
Hom_{\f}(B_1,B_2)=\un_{x\in\set}(x,Hom_{\f}(E_x^1,E_x^2))
\end{equation*}
We have a $G$-action on the vectorbundle $Hom_{\f}(B_1,B_2)$ given by \linebreak $g(x,\phi_x)=(gx,g(\phi_x))$
where we define
\begin{equation*}
g(\phi_x)([e_1]_{gx})=g(\phi_x(g^{-1}[e_1]_{gx}))
\end{equation*}
This induces the structure of a left A-module on the set $\gm(Hom_{\f}(B_1,B_2))$ in the usual
way.

\begin{proposition}
$Hom_k(\e_1,\e_2)\iso\gm(Hom_{\f}(B_1,B_2))$
\end{proposition}
\begin{proof}
Define a map 
\begin{align*}
F:Hom_k(\e_1,\e_2)&\ft\gm(Hom_{\f}(B_1,B_2)),\\
\phi&\mt F\phi
\end{align*}
as $(F\phi)(x)(v_x^1)=\phi(s)(x)$ where $s\in\gm(B_1)$ satisfies $s(x)=v_x^1$. This is well defined 
because if $s,s'\in\gm(B_1)$ and $s(x)=s'(x)$ for the given $x\in\set$ then
\begin{align*}
\phi(s)(x)&=(\delta_x\phi(s))(x)=(\phi(\delta_xs))(x).\\
\phi(s')(x)&=(\delta_x\phi(s'))(x)=(\phi(\delta_xs'))(x).
\end{align*}
But
\begin{align*}
(\delta_xs)(y)&=\delta_x(y)s(y)=\delta_{xy}v_x^1,\\
(\delta_xs')(y)&=\delta_x(y)s'(y)=\delta_{xy}v_x^1.
\end{align*}
Therefore $\delta_xs=\delta_xs'$ and we have $\phi(s)(x)=(\phi(\delta_xs))(x)=(\phi(\delta_xs'))(x)
=\phi(s')(x)$. So $F$ is well defined.
Assume that $F\phi_1=F\phi_2$. Let $s\in\gm(B_1)$, then $\phi_1(s)(x)=(F\phi_1)(x)(s(x))=
(F\phi_2)(x)(s(x))=\phi_2(s)(x)$ for all $x\in\set$. But then $\phi_1=\phi_2$ and $F$ is
injective.
Let $\gamma\in\gm(Hom_{\f}(B_1,B_2))$ be given. Define a map $\phi:\gm(B_1)\ft\gm(B_2)$ by
$\phi(s)(x)=\gamma(x)(s(x))$. Then we have 
\begin{align*}
\phi(s+s')(x)&=\gamma(x)(s(x)+s'(x))\\
&=\gamma(x)(s(x))+\gamma(x)(s'(x))=\phi(s)(x)+\phi(s')(x),\\
\phi(fs)(x)&=\gamma(x)(f(x)s(x))=f(x)\gamma(x)(s(x))\\
&=f(x)\phi(s)(x)=(f\phi(s))(x).
\end{align*}
So we have that $\phi\in Hom_k(\gm(B_1),\gm(B_2))$ and also
\begin{align*}
F(\phi)(x)(v_x^1)&=\phi(s)(x)=\gamma(x)(s(x))=\gamma(x)(v_x^1).
\end{align*}
Therefore we have that $F\phi=\gamma$ and $F$ is surjective. Furthermore we have
\begin{align*}
F(f\phi)(x)(v_x^1)&=(f\phi)(s)(x)=(f(\phi(s)))(x)\\
&=f(x)(\phi(s)(x))=f(x)(F(\phi)(x)(v_x^1))=(f(x)F(\phi)(x))(v_x^1)=(fF(\phi))(x)(v_x^1),\\
F(g\phi)(x)(v_x^1)&=(g\phi)(s)(x)=(g(\phi(g^{-1}s)))(x)=g(\phi(g^{-1}s)(g^{-1}x))\\
&=g(F(\phi)(g^{-1}x)(g^{-1}v_x^1))=(g(F(\phi)(g^{-1}x)))(v_x^1)\\
&=(gF(\phi))(x)(v_x^1).
\end{align*}
So $F$ is a A-module morphism.
\end{proof}

Let $\e\iso\gm(B)$ be a given GF-difference equation  where $B=\un_{x\in\set}(x,E_x)$ is a vectorbundle.
Define the dual vector bundle $B^*=\un_{x\in\set}(x,E_x^*)$. Then as a special case of the
previous proposition we have.

\begin{proposition}
$\e^*\iso\gm(B^*)$.
\end{proposition}

\subsection{The geometric description of A-morphisms}

Let $\gm(B),\gm(B')$ be two GF-difference equations with corresponding vectorbundles $B=\un_{x\in\set}(x,E_x)
,B'=\un_{x\in\set}(x,E_x')$ and let $\phi\in Hom_A(\gm(B),\gm(B'))$. Define
a map $F_{\phi}:B\ft B'$ by 
\begin{equation*}
F_{\phi}(x,v_x)=(x,F_{\phi}^x)
\end{equation*}

where $F_{\phi}^x(v_x)=\phi(s)(x)$ and $s\in\gm(B)$ is any section satisfying $s(x)=v_x$.

\begin{proposition}
$F_{\phi}^x:E_x\ft E'_x$ is well defined.
\end{proposition}
\begin{proof}
Assume $s(x)=s'(x)=v_x$. then
\begin{align*}
\phi(s)(x)&=(\delta_x\phi(s))(x)=\phi(\delta_xs)(x),\\
\phi(s')(x)&=(\delta_x\phi(s'))(x)=\phi(\delta_xs')(x),
\end{align*}
and $(\delta_xs)(y)=\delta_x(y)s(y)=\delta_x(y)s'(y)=(\delta_xs')(y)$ for all
$y\in\set$. This means that $\delta_xs=\delta_xs'$ and so $\phi(\delta_xs)=\phi(\delta_xs')$
and we can conclude that $\phi(s)(x)=\phi(s')(x)$.
\end{proof}

\begin{proposition}
$F_{\phi}^x$ is $\f$-linear.
\end{proposition}
\begin{proof}
Let $v_x,u_x\in E_x$ and let $s,t\in\gm(B)$ be any sections such that $s(x)=v_x,t(x)=u_x$.
Then $(s+t)(x)=v_x+u_x$ and we have
\begin{align*}
F_{\phi}^x(v_x+u_x)&=\phi(s+t)(x)=\phi(s)(x)+\phi(t)(x)=F_{\phi}^x(v_x)+F_{\phi}^x(u_x).
\end{align*}
Let $a\in\f$ and $v_x\in E_x$. Let $s\in\gm(B)$ be any section such that $s(x)=v_x$. Then
$(as)(x)=a(s(x))=av_x$ and we have
\begin{align*}
F_{\phi}^x(av_x)&=\phi(as)(x)=(a\phi(s))(x)=a(\phi(s)(x))=aF_{\phi}^x(v_x),
\end{align*}
\end{proof}

\begin{proposition}
$F_{\phi}^y\circ g=g\circ F_{\phi}^{g^{-1}y}$ for all $g\in G$ and $y\in\set$.
\end{proposition}
\begin{proof}
Let $y\in\set,g\in G$ and $v_{g^{-1}y}\in E_{g^{-1}y}$. Let $s\in\gm(B)$ be any section such that
$s(g^{-1}y)=v_{g^{-1}y}$. Then we have
\begin{align*}
F_{\phi}^y(gv_{g^{-1}y})&=F_{\phi}^y(gs(g^{-1}y))=F_{\phi}^y((gs)(y))=\phi(gs)(y)\\
&=(g\phi(s))(y)=g(\phi(s)(g^{-1}y))=g(F_{\phi}^{g^{-1}y}(s(g^{-1}y)))=g(F_{\phi}^{g^{-1}y}(v_{g^{-1}y}).
\end{align*}
\end{proof}

The previous three propositions show that a morphism of GF-difference equations is a family
of $\f$-linear maps that are related at different points as described in the last proposition.
Let $H_y$ be the isotropy group of the point $y\in\set$. As a special case of the last proposition we have

\begin{corollary}
$F_{\phi}^y\circ h=h\circ F_{\phi}^y$ for all $h\in H_y$ and $y\in\set$.
\end{corollary}

So the maps $F_{\phi}^y$ commutes with the action of the isotropy group at each point
and are $\f[H_x]$-module morphisms on the fibre above the point.
Properties of the morphisms $\phi\in Hom(\gm(B),\gm(B'))$ is transferred to the family
of maps $F_{\phi}^x$.

\begin{proposition}
Let $x\in\set$ be some point in $\set$. Then $\phi\in Hom_A(\gm(B),\gm(B'))$ is surjective
if and only if $F_{\phi}^x:E_x\ft E_x'$ is surjective.
\end{proposition}
\begin{proof}
Assume $\phi$ is surjective. Let $v_x'\in E_x'$ be given. Then there exists $\gamma\in\gm(B')$ such
that $\gamma(x)=v_x'$. Let $s\in\gm(B)$ be such that $\phi(s)=\gamma$. Let $v_x=s(x)$. Then
\begin{align*}
F_{\phi}^x(v_x)=\phi(s)(x)=\gamma(x)=v_x'
\end{align*}
so $F_{\phi}^x$ is surjective.

Assume that $F_{\phi}^x$ is surjective. Let $y\in \set$ and let $v_y'\in E_y'$. There exists $g\in G$ such
that $gx=y$. Define $v_x'=g^{-1}v_y'\in E_x'$. Then there exists $v_x\in E_x$ such that 
$F_{\phi}^x(v_x)=v_x'$. Define $v_y=gv_x$. Then we have
\begin{align*}
F_{\phi}^y(v_y)=F_{\phi}^y(gv_x)=gF_{\phi}^x(v_x)=gv_x'=v_y'
\end{align*}
so $F_{\phi}^y$ is surjective for all $y\in\set$. Let $\gamma\in\gm(B')$. Then $\gamma(y)=v_y\in E_y'$ for
 all y. For each y there then exists $v_y\in E_y$ such that $F_{\phi}^y(v_y)=v_y'$. Define 
$s\in\gm(B)$ by $s(y)=v_y$. Then we have
\begin{align*}
\phi(s)(y)=F_{\phi}^y(v_y)=v_y'=\gamma(y)
\end{align*}
so $\phi(s)=\gamma$ and $\phi$ is surjective.
\end{proof}

\begin{proposition}
Let $x\in\set$ be any point in $\set$. Then $\phi\in Hom_A(\gm(B),\gm(B'))$ is injective if
and only if $F_{\phi}^x:E_x\ft E_x'$ is injective.
\end{proposition}
\begin{proof}
Assume that $\phi\in Hom_A(\gm(B),\gm(B'))$ is not injective. Then there exists $s\in\gm(B)$
,$s\ne 0$ such that $\phi(s)=0$. There is at least one point $y\in\set$ such that $s(y)=v_y\ne 0$. Then
\begin{align*}
F_{\phi}^y(v_y)=\phi(s)(y)=0
\end{align*}
so $F_{\phi}^y$ is not injective. Let $g\in G$ be such that $gx=y$. Let $g:E_x\ft E_y$ be the corresponding
invertible fibermap. Define $v_x=g^{-1}v_y$. Then $v_x\ne 0$,$v_x\in E_x$ and
\begin{align*}
F_{\phi}^x(v_x)=F_{\phi}^x(g^{-1}v_y)=g^{-1}F_{\phi}^y(v_y)=0
\end{align*}
so $F_{\phi}^x$ is not injective.

Assume that $\phi$ is injective. Let $v_x\in E_x$ and assume that $F_{\phi}^x(v_x)=0$. Let $s\in\gm(B)$ be
any section such that $s(x)=v_x$. Define $\gamma\in\gm(B)$ by $\gamma=\delta_x s$. Then
\begin{align*}
\phi(\gamma)(y)=\phi(\delta_x s)(y)=\delta_x(y)\phi(s)(y)=\delta_{xy}F_{\phi}^x(v_x)=0
\end{align*}

for all $y\in\set$. But then $\gamma=0$ and so $v_x=s(x)=\gamma(x)=0$ and we conclude that $F_{\phi}^x$
is injective.
\end{proof}

Combining the previous  propositions we have

\begin{corollary}\label{iso1}
Let $x\in\set$ be any point in $\set$. Then $\phi\in Hom_A(\gm(B),\gm(B'))$ is an
A-module isomorphism if and only if $F_{\phi}^x:E_x\ft E_x'$ is a $\f[H]$-module
isomorphism.
\end{corollary}

Any A-morphism gives us a family of $\f$-linear maps with the properties described. Any such
family will in fact come from a A-morphism of modules.

For each $x\in\set$ let $F^x:E_x\ft E_x'$ be a $\f$-linear map. Assume that the members of the family 
are related through

\begin{align*}
F_{\phi}^x\circ g=g\circ F_{\phi}^{g^{-1}x}
\end{align*}

for all $g\in G$ and $x\in\set$.

 Define a map $\phi:\gm(B)\ft\gm(B')$ by
\begin{align*}
\phi(s)(x)=F^x(s(x)).
\end{align*}

\begin{proposition}\label{hom1}
$\phi\in Hom_A(\gm(B),\gm(B'))$
\end{proposition}
\begin{proof}
Let $s,t\in\gm(B)$. Then
\begin{align*}
\phi(s+t)(x)&=F^x((s+t)(x))=F^x(s(x)+t(x))\\
&=F^x(s(x))+F^x(t(x))=\phi(s)(x)+\phi(t)(x).
\end{align*}
Let $s\in\gm(B)$ and $f\in k$. Then
\begin{align*}
\phi(fs)(x)&=F^x((fs)(x))=F^x(f(x)s(x))\\
&=f(x)F^x(s(x))=f(x)\phi(s)(x)=(f\phi(s))(x).
\end{align*}
Let $s\in\gm(B)$ and $g\in G$. Then
\begin{align*}
\phi(gs)(x)&=F^x((gs)(x))=F^x(g(s(g^{-1}x)))\\
&=g(F^{g^{-1}x}(s(g^{-1}x)))=g(\phi(s)(g^{-1}x))
=(g\phi(s))(x)
\end{align*}
\end{proof}

In general a submodule of a finitely generated module does not have to be finitely
generated. We will now show that for the category $\GFE$ all submodules are
in fact GF-difference equations.

\begin{proposition}
Let $\e$ and $\e'$ be GF-difference equations and let $\phi\in Hom_A(\e,\e')$ be
a A-module morphism. The $Im\phi$ and $Ker\phi$ are GF-difference equations.
\end{proposition}
\begin{proof}
We know that $im\phi$ and $ker\phi$ are submodules of geometric modules and are therefore geometric.
Let $\{e_i\}_{i=1}^{n}$ be a set of generators for $\e$. Then $\{\phi(e_i)\}_{i=1}^n$
is a finite set of generators for $im\phi$. So $im\phi$ is a GF-difference equation.
We know that $\e\iso\gm(B)$ for some vectorbundle $B=\un_{y\in\set}(y,E_y)$. We know that
$ker\phi$ is a geometric submodule so we have an injective A-module morphism 
$T:ker\phi\hookrightarrow\gm(B')\sub\gm(B)$ where $B'=\un_{y\in\set}(y,V_y)$ is the
subbundle with fibers
\begin{align*}
V_y=\{s(y)\mid s\in ker\phi\}.
\end{align*}
Let $\{F_{\phi}^y\}_{y\in\set}$ be the family of maps corresponding to $\phi$ and let $v_y\in V_y$.
Then $v_y=s(y)$ for some $s\in ker\phi$ and we have $F_{\phi}^y(v_y)=\phi(s)(y)=0$. Let $\gamma\in\gm(B')$.
Then $\gamma(y)\in V_y$ for all $y\in\set$ so that $\phi(\gamma)(y)=F_{\phi}^y(\gamma(y))=0$ and as a
concequence $\gamma\in ker\phi$ and the map T is surjective. But then we have proved that $ker\phi\iso\gm(B')$ and $ker\phi$ is finitely generated and therefore a GF-difference equation.
\end{proof}
 
\begin{corollary}
Let $\e$ be a GF-difference equation and $\e'\sub\e$ a submodule. Then $\e'$ is a GF-difference equation.
\end{corollary}
\begin{proof}
We know that $\e$ and $\e/\e'$ are GF-difference equations. Let $\phi:\e\ft\e/\e'$ be the natural
projection. Then $\phi$ is a A-module morphism and $\e'=ker\phi$.
\end{proof}

\subsection{The general structure of $\GFE$}

We are now ready to give a characterization of the structure of the category $\GFE$. All the
structural properties follow from the following proposition.
\begin{proposition}
Let $\e$ be a GF-difference equation. Then $\e$ is both Artinian and Noetherian.
\end{proposition}
\begin{proof}
Let $\e_1\sub\e_2\sub\e_3\cdots$ be a ascending chain of submodules in $\e$. Then in particular this
is a chain of free k-modules. But $\e$ has finite dimension over k so the chain must stop and the
module $\e$ is Noetherian.
Similarly let $\cdots\sub\e_3\sub\e_2\sub\e_1$ be a decending chain of submodules of $\e$. Then in particular
it is a decending chain of free k-modules. But the dimension of any module over k is nonnegative so the chain
must stop.
\end{proof}

A GF-difference equation is $\emph{simple}$ if it contains no GF-difference equations as submodules and 
$\emph{indecomposable}$ if
 it can not be written as a direct sum of GF-difference equations. A composition series for a GF-difference 
equation is a finite filtering $0\sub\e_1\sub\e_2\cdots\e$ of the equation $\e$ such that the composition
factors $\e_i/\e_{i-1}$ are simple equations. Because of the previous proposition and the general structure
 theory for modules \cite{row} we have
\begin{theorem}
Let $\e$ be a GF-difference equation. Then $\e$ has a composition series and all composition series for $\e$ 
has the same number of elements in the filtration and the composition factors are the same up to isomorphism.
Any GF-difference equation can be written as a direct sum of a finite number of indecomposable equations.
\end{theorem}

This theorem is a combination of the Jordan-H\o\/lder theorem and the Krull-Smidth teorem. This Theorem
reduce the problem of solving GF-difference equations to the study of indecomposable equations. Furthermore
it show that all indecomposable equations  are related to simple equations through a finite set of
simple decomposition factors. The first problem then is to find the simple equations and the next is
to construct the indecomposable equations using a set of composition factors. This last
problem is essentially the problem of quantization. Even finding the simple equations is in general
not a trivial task in the category $\GFE$. We will however now proceed to prove a Theorem that show that the Category $\GFE$ is
equivalent to a category where the problem of finding simple and indecomposable objects is more
approachable.

\section{The equivalence Theorem}

We will first construct a special class of GF-difference equations and then show that all
GF-difference equations are in fact of this type. Let $x\in\set$ and let $H_x$ be the isotropy group of that
point. All such groups for different points $x$ are isomorphic. We will usually supress the
point we are referring to and just write $H=H_x$. For each $y\in\set$ define the set $yH$ by
\begin{equation*}
yH=\{g\in G\mid gx=y\}.
\end{equation*}
We evidently have $gH=yH$ where $g\in G$ satisfy $gx=y$ so the sets $yH$ are just the
left cosets of $H$ in $G$. Let $V$ be a finite dimensional $\f[H]$-module and form the
trivial bundle $\set\times V$
\begin{equation*}
\set\times V=\un_{y\in\set}(y,V),
\end{equation*}
Let $\sigma$ be a transversal to the partitioning of $G$ by the classes $yH$
\begin{equation*}
\sigma(y)\in yH\quad\text{$\forall y\in\set$},
\end{equation*}
Let $\gm(G\times V)$ be the k-module of sections in the trivial bundle $\set\times V$. We will now 
define an action of $G$ in this module of sections as
\begin{equation*}
g(y,v)=(gy,\sigma(gy)^{-1}g\sigma(y)v).
\end{equation*}
This gives a $G$-action.
\begin{proposition}
$(gg')(y,v)=g(g'(y,v))$
\end{proposition}
\begin{proof}
\begin{align*}
(gg')(y,v)&=((gg')y,\sigma((gg')y)(gg')\sigma(y)v)\\
&=(g(g'y),\sigma(g(g'y))g\sigma(g'y)\sigma(g'y)^{-1}g'\sigma(y)v)\\
&=g(g'y,\sigma(g'y)g'\sigma(y)v)=g(g'(y,v)).
\end{align*}
\end{proof}
We let this action in the bundle induce an action on the module of sections in the bundle
in the usual way.
\begin{equation*}
(gs)(y)=g(s(g^{-1}y))
\end{equation*}
This gives the set of section in the bundle $\set\times V$ the structure of a GF-difference equation.
It appears as if each choise of transversal $\sigma$ gives a different action
on the bundle and so a different A-module structure on the set $\gm(\set\times V)$. They are
however all isomorphic. Let $\gm(\set\times V)$ and $\gm(\set\times V)'$ be the modules corresponding
to the choise of two transversals $\sigma$ and $\sigma'$. Then we have
\begin{proposition}
$\gm(\set\times V)\iso\gm(\set\times V)'$
\end{proposition}
\begin{proof}
For each $y\in\set$ there exists a $\gamma(y)\in H$ such that $\sigma(y)=\sigma'(y)\gamma(y)$.
Define a map $\phi_{\gamma}:\gm(\set\times V)\ft\gm(\set\times V)'$ by $\phi_{\gamma}(s)(y)=
\gamma(y)^{-1}(s(y))$ This map is clearly bijective with inverse $\phi_{\gamma^{-1}}$ where
$\phi_{\gamma^{-1}}(s)(y)=\gamma(y)(s(y))$. We also have
\begin{align*}
\phi_{\gamma}(fs)(y)&=\gamma(y)^{-1}((fs)(y))=\gamma(y)^{-1}(f(y)s(y))=f(y)(\gamma(y)^{-1}(s(y)))\\
&=f(y)\phi_{\gamma}(s)(y)=(f\phi_{\gamma}(s))(y),\\
\phi_{\gamma}(gs)(y)&=\gamma(y)^{-1}((gs)(y))=\gamma(y)^{-1}(g(s(g^{-1}y)))\\
&=\gamma(y)^{-1}(\sigma(y)^{-1}g\sigma(g^{-1}y)s(g^{-1}y))\\
&=\gamma(y)^{-1}\sigma(y)^{-1}g\sigma(g^{-1}y)\gamma(g^{-1}y)\gamma(g^{-1}y)^{-1}s(g^{-1}y)\\
&=(\sigma(y)\gamma(y))^{-1}g(\sigma(g^{-1}y)\gamma(g^{-1}y))\gamma(g^{-1}y)^{-1}s(g^{-1}y)\\
&=\sigma'(y)^{-1}g\sigma'(g^{-1}y)(\phi_{\gamma}(s))(g^{-1}y)\\
&=(g\phi_{\gamma}(s))(y).
\end{align*}
So the map $\phi_{\gamma}$ is also a A-module morphism and the proof is complete.
\end{proof}

The constructed class of GF-difference equations in fact includes all GF-difference equations.

\begin{theorem}\label{Theorem}
Let $\e$ be any GF-difference equation. Then $\e\iso\gm(\set\times V)$ for some $\f[H]$-module of 
finite dimension over $\f$.
\end{theorem}
\begin{proof}
We know that $\e\iso\gm(B)$ for some vectorbundle $B$
\begin{equation*}
B=\un_{y\in\set}(y,E_y).
\end{equation*}
Let $\sigma$ be a transversal to the classes $yH$, Then $\sigma(y)x=y$ and so the action of $G$
on the bundle gives us the lift $\sigma(y):E_x\ft E_y$ and this map is an isomorphism since $G$
is a group. Define $V=E_x$. Then $V$ is an finite dimensional $H$ space. Let $\phi$ be the map
\begin{align*}
\phi:B&\ft\set\times V\\
(y,v_y)&\mt(y,\sigma(y)^{-1}v_y).
\end{align*}
This map is clearly an isomorphism of vectorbundles and it commutes with the $G$ action
\begin{align*}
\phi(g(y,v_y))&=\phi(gy,gv_y)=(gy,\sigma(gy)^{-1}gv_y)\\
&=(gy,\sigma(gy)^{-1}g\sigma(y)\sigma(y)^{-1}v_y)=g\phi(y,v_y).
\end{align*}
Define a map $F_{\phi}$ by
\begin{align*}
F_{\phi}:\gm(B)&\ft\gm(\set\times V)\\
F_{\phi}(s)(y)&\mt\phi(s(y))
\end{align*}
$F_{\phi}$ is clearly bijective with inverse $F_{\phi^{-1}}$ and we also have
\begin{align*}
F_{\phi}(fs)(y)&=\phi((fs)(y))=\phi(f(y)s(y))=f(y)\phi(s(y))\\
&=f(y)(F_{\phi}(s)(y))=(fF_{\phi}(s))(y),\\
F_{\phi}(gs)(y)&=\phi((gs)(y))=\phi(g(s(g^{-1}y)))\\
&=g(\phi(s(g^{-1}y)))=g(F_{\phi}(s)(g^{-1}y))=(gF_{\phi}(s))(y).
\end{align*}
So $F_{\phi}$ is a A-module morphism and the proof is complete.
\end{proof}
 
Let $\f[H]-finmod$ be the category of modules over $\f[H]$ with finite dimension
over $\f$ with direct sum and tensor product and dual over $\f$ defined as is usual in
representation theory. For this category all main structural theorems for the decomposition
of modules apply so that all such modules have compositionseries and can be written as a finite 
direct sum of indecomposables. We will now proceed to show that the category $\GFE$ and $\f[H]-finmod$
are in fact equivalent as categories. From a structural point of view we will not distinguish between 
isomorphic objects 
and will therefore prove the equivalence by showing the isomorphism of the Grotendick
algebra \cite{grot} $A_G$ for $\GFE$ and $A_H$ for $\f[H]-finmod$. The algebra structure in $A_G$ and $A_H$ is
the one induced from direct sum and tensor product in the corresponding categories. In addition to
the usual algebra structure we have a conjugation map induced from the dual in the categories.
Define a map on $T:A_H\ft A_G$ by

\begin{equation*}
T([V])=[\gm(\set\times V)]
\end{equation*}

where elements in the Grotendick algebras are denoted by square brackets of elements in the
corresponding categories.

\begin{proposition}
T is well defined
\end{proposition}
\begin{proof}
Assume $[V]=[U]$. Then there exists a $\f[H]$-module isomorphism $F^x:V\ft U$. So
$F^x$ is an $\f$-isomorphism and $F^x(hv)=hF^x(v)$ for all $h\in H$. Denote the fiber over
 $y\in\set$  of the vectorbundles $B_1=\set\times V$ and $B_2=\set\times U$ by $V_y$ and $U_y$.
For each $y\in\set$ define a map $F^y:V_y\ft U_y$ by $F^y(v_y)=g(F^x(g^{-1}v_y))\in U_y$ where g is any 
element in G such that $gx=y$. The family of maps $\{F^y\}_{y\in\set}$ is well defined 
because if $g_1$ is another element in G such that $g_1x=y$ then $g_1=gh$ for some $h\in H$ and we have
\begin{align*}
g_1(F^x(g_1^{-1}v_y))&=(gh)(F^x((gh)^{-1}v_y))=g(h(F^x(h^{-1}g^{-1}v_y)))=g(F^x(g^{-1}v_y)).
\end{align*}
The family $\{F^y\}_{y\in \set}$ satisfies all requirements in proposition (\ref{hom1}) and
therefore determines a A-morphism, $\phi:\gm(\set\times V)\ft\gm(\set\times\ U)$.
 From the construction
we observe that each member of the family $\{F^y\}_{y\in\set}$ is a $\f[H]$-module isomorphism.
We therefore can conclude that the map $\phi$ is an A-isomorphism so that $[\gm(\set\times V)]
=[\gm(\set\times U)]$.
\end{proof}

\begin{proposition}
$T:A_H\ft A_G$ is a bijection.
\end{proposition}
\begin{proof}
Let $[\e]$ be any element in $A_G$. From theorem \ref{Theorem} we know that there exists
a $V$ in $\f[H]-finmod$ such that $[\e]=[\gm(\set\times V)]$. The $T([V])=[\e]$ so that
$T$ is surjective. Assume that $T([V])=T([U])$. This means that there exists a A-module
isomorphism $\phi:\gm(\set\times V)\ft\gm(\set\times U)$. Let $y\in\set$ be any point.Then 
proposition\ref{iso1} show that there exists a $\f[H]$-module isomorphism $F_{\phi}^y:V\ft U$. But then 
$[V]=[U]$ and T is injective.
\end{proof}

Rewriting some of the results proved earlier we find that T is a structure preserving map.

\begin{proposition}
The map T is structure preserving
\begin{align*}
T([U]+[V])&=T([U])+T([V])\\
T([U][V])&=T([U])T([V])\\
T([U]^*)&=(T([U])^*
\end{align*}
\end{proposition}

This relation between the categories $\GFE$ and $\f[H]-finmod$ gives us a way to find all
 indecomposable equations of
finite type in $\GFE$ from
the indecomposable $\f[H]$-modules of finite dimension over $\f$.
\begin{proposition}
$\e$ is a indecomposable GF-difference equation of finite type if and only if 
$\e\iso\gm(\set\times V)$ where is a indecomposable $\f[H]$-module of finite dimension over $\f$.
\end{proposition}
\begin{proof}
Let $\e$ be indecomposable. We know that $\e\iso\gm(\set\times V)$ for some $\f[H]$-module V.
Assume that V is decomposable so that $V\iso V_1\ds V_2$. Define $\e_1=\gm(\set\times V_1)
,\e_2=\gm(\set\times V_2)$. Then 
\begin{align*}
[\e]&=T([V])=T([V_1\ds V_2])=T([V_1]+[V_2])\\
&=T([V_1])+T([V_2])=[\e_1]+[\e_2]=[\e_1\ds\e_2].
\end{align*}
So that $\e\iso\e_1\ds\e_2$ and $\e$ is decomposable. This is a contradiction so that V
 is indecomposable. 
Let V be indecomposable and of finite dimension over $\f$. Define $\e=\gm(\set\times V)$.
Then $\e$ is a GF-difference equation of finite type and $T([V])=[\e]$. Assume that $\e$ is decomposable.
Then $\e\iso\e_1\ds\e_2$. We know that there exists $\f[H]$-modules $V_1,V_2$ of finite dimension
over $\f$ such that $\e_1\iso\gm(\set\times V_1),\e_2\iso\gm(\set\times V_2)$. Then we have
\begin{align*}
T[V]&=[\e]=[\e_1\ds\e_2]=[\e_1]+[\e_2]=T([V_1])+T([V_2])\\
&=T([V_1]+[V_2])=T([V_1\ds V_2]).
\end{align*}
But T is injective so that $[V]=[V_1\ds V_2]$. Then $V\iso V_1\ds V_2$ and V is decomposable.
This is a contradiction.
\end{proof}

For simple equations of finite type we have
\begin{proposition}
$\e$ is a simple GF-difference equation of finite type if and only if $\e\iso\gm(\set\times V)$
where V is a simple $\f[H]$-module.
\end{proposition}
\begin{proof}
Let $\e$ be a simple GF-difference equation of finite type. We know that $\e\iso\gm(\set\times V)$ for
some $\f[H]$-module V of finite dimension over $\f$. Assume that V has a submodule $V'$. let
$\e'=\gm(\set\times V')$. Then $\e'$ is a submodule of $\e$ so $\e$ is not simple. This is a contradiction.
Assume that V is a simple $\f[H]$-module. Let $\e=\gm(\set\times V)$ and assume that 
$\e$ is not simple so that it has a submodule $\e'$. But then $\e'\iso\gm(\set\times V')$ where $V'$ is
a submodule of V. This is a contradiction.
\end{proof}

The simple $\f[H]$-modules are in general not easy to find. For the case when the 
isotropy group is finite and the character of the field does not divide the order of
the group H, the algebra $\f[H]$ is semisimple  and the full power of the theory
of characters \cite{isa} can be applied. Even in the case when the character of $\f$ does divide
the order of the group, the modular case, powerful tools are available.

\section{The projection formula for GF-difference equations}

The Frobenius projection formula \cite{frob} can be generalized to apply to GF-difference equations
when the isotropy group is finite and the underlying field is $\zc$. This formula greatly 
simplifies the solution process
when it applies.
Let $\e\iso\gm(B)$ be any GF-difference equation where $B=\un_{y\in\set}(y,E_y)$ is a
vectorbundle over $\set$. We know that the fiber over $y\in\set$ is a $\f[H_y]$-module.
Denote the character of this module by $\chi_y$. There is a relation between characters
at different points.
\begin{proposition}
$\chi_y(h_y)=\chi_{gy}(gh_yg^{-1})$ for all $y\in\set,g\in G$ and $h_y\in H_y$.
\end{proposition}
\begin{proof}
Let $h_y\in H_y$ and let $g\in G$. Then $h_y=g^{-1}gh_yg^{-1}g=g^{-1}(gh_yg^{-1})g$
so we have
\begin{align*}
\chi_y(h_y)=tr(h_y)=tr(g^{-1}(gh_yg^{-1})g)=tr(gh_yg^{-1})=\chi_{gy}(gh_yg^{-1}).
\end{align*}
\end{proof}

Assume that $\e'\iso\gm(B')$ is a submodule of $\e\iso\gm(B)$. Then $B'\sub B$ so
that $E'_y\sub E_y$ as a $\f$-linear subspace for all $y\in\set$. For each $y\in\set$
let $\Pi_{\e'}(y):E_y\ft E_y$ be the Frobenius map
\begin{align*}
\Pi_{\e'}(y)v_y=\frac{dim E'_y}{\mid H_y\mid}\sum_{h_y\in H_y}\chi_y(h_y^{-1})h_yv_y. 
\end{align*}

Here $\chi_y$ is the character of the $\f[H]$-module $E_y'$.

Each $h_y$ is a $\f$-linear map so clearly $\Pi_{\e'}(y)$ is a $\f$-linear map for each
$y\in\set$. We also have

\begin{proposition}
$\Pi_{\e'}(y)(gv_{g^{-1}y})=g\Pi_{\e'}(g^{-1}y)(v_{g^{-1}y})$ for all $y\in\set$ and $g\in G$.
\end{proposition}
\begin{proof}
Let $c_y=dim E_y/\mid H_y\mid$. Then $c_y=c_{y'}$ for all $y,y'\in\set$ since $B=\un_{y\in\set}(y,E_y)$
is a vectorbundle and $H_y\iso H_{y'}$ for all $y,y'\in\set$.
We have
\begin{align*}
\Pi_{\e'}(y)(g^{-1}v_y)&=c_y\sum_{h_y\in H_y}\chi_y(h_y^{-1})h_y(gv_{g^{-1}y})=c_y\sum_{h_y\in H_y}\chi_y(h_y^{-1})
g(g^{-1}h_yg)v_{g^{-1}y})\\
&=c_yg\sum_{h_y\in H_y}\chi_y(h_y^{-1})(g^{-1}h_yg)v_{g^{-1}y}\\
&=c_yg\sum_{h_{g^{-1}y}\in H_{g^{-1}y}}
\chi_y(gh_y^{-1}g^{-1})h_{g^{-1}y}v_{g^{-1}y}\\
&=c_{g^{-1}y}g\sum_{h_{g^{-1}y}\in H_{g^{-1}y}}
\chi_{g^{-1}y}(h_{g^{-1}y}^{-1})h_{g^{-1}y}v_{g^{-1}y}\\
&=\Pi_{\e'}(g^{-1}y)(v_{g^{-1}y}).
\end{align*}
\end{proof}

Define a map $\Pi_{\e'}$ on $\e\iso\gm(B)$ by
\begin{align*}
\Pi_{\e'}(e)(y)=\Pi_{\e'}(y)(e(y)).
\end{align*}
Then we can conclude from the previous proposition that
\begin{corollary}
$\Pi_{\e'}\in Hom_A(\e,\e)$.
\end{corollary}

In general $\Pi_{\e'}$ is not a projection on $\e'$. Assume that $\f=\zc$. Then $\e=\sum_in_i\s_i$
where all $\s_i$ are simple A-modules. Then we have
\begin{proposition}
$\Pi_{\s_i}\mid_{\s_j}=\delta_{ij}id$
\end{proposition}
\begin{proof}
$\Pi_{\s_i}\mid_{\s_j}\in Hom_A(\s_j,\s_j)$. But this means that $\Pi_{\s_i}\mid_{\s_j}(y)
\in Hom_{\zc}(\s_{jy},\s_{jy})$ for all $y\in\set$. From the Schur lemma we can conclude that
$\Pi_{\s_i}\mid_{\s_j}(y)=\lambda_{ij}(y)id_y$. But then
\begin{align*}
dim \s_{jy} \lambda_{ij}(y)=Tr(\Pi_{\s_i}\mid_{\s_j}(y))=\frac{dim \s_{iy}}{\mid H_y\mid}
\sum_{h_y\in H_y}\chi_{iy}(h_y^{-1})\chi_{jy}(h_y)=\delta_{ij}dim \s_{iy}
\end{align*}
So we have $\lambda_{ij}(y)=\delta_{ij}$ and the proof is complete.
\end{proof}

\begin{corollary}
$\Pi_{\s_i}$ is a projection of $\e$ onto $n_i\s_i\sub\e$.
\end{corollary}

The projection formula can be used to simplify the solution process for GF-difference
equations.
\begin{proposition}
$\phi$ is a solution of $\e$ of type $\s$ if and only of $\phi=\Pi_{\s}^*(\psi)$
for some solution $\psi$  of $\Pi_{\s}(\e)$ of type $\s$.
\end{proposition}
\begin{proof}
Let $\psi\in Hom_A(n_i\s_i,\s_i)$. Then $\Pi_{\s_i}^*(\psi)=\psi\circ\Pi_{\s_i}\in Hom_A(\e,\s_i)$
and so $\Pi_{\s_i}^*(\psi)$ is a solution of $\e$ of type $\s_i$.
Conversely let $\phi\in Hom_A(\e,\s_i)$ be a solution of $\e$ of type $\s_i$. We know by the
structure theorem that $\e\iso\sum_jn_j\s_j$. By the Schur lemma $\phi\mid_{n_j\s_j}=0$ for
$i\ne j$. Let $\psi=\phi\mid_{n_i\s_i}$, then $\psi\in Hom_A(n_i\s_i,\s_i)$ is a solution of
$n_i\s_i$ of type $\s_i$ and for any $e=\sum_je_j$ in $\e\iso\sum_jn_j\s_j$ we have
\begin{align*}
\Pi_{\s_i}^*(\psi)(e)&=\psi(\Pi_{\s_i}(e))=\psi(e_i)=\phi(e_i)=\phi(e).
\end{align*}
\end{proof}

\section{Coordinate description of GF-difference equations}

Let $\e$ be any GF-difference equation. The structure of $\e$ is essentially determined by the action
of $G$ on the k-module $\e$. Let $\{e_i\}_{i=1}^n$ be a k-basis for $\e$. Define a set of matrices
$\e^g\in Mat(n,k)$ by
\begin{equation*}
ge_i=\sum_j\e^g_{ij}e_j
\end{equation*}
The set of matrices $\{\e^g\}_{g\in G}$ determines the $G$-action on $\e$ with respect to the 
gives k-basis. They formally play the same role as the connection symbols in differential
geometry and we will call them $\emph{the connection of the given GF-difference equation}$.
In general $\e^{gg'}\ne\e^g\e^{g'}$ so the relation $g\ft\e^g$ is not a representation of $G$.
We have however the following result.

\begin{proposition}
 $\e^{gg'}=g(\e^{g'})\e^{g}$.
\end{proposition}
\begin{proof}
\begin{align*}
(gg')e_i&=g(g'e_i)=g\sum_j\e^{g'}_{ij}e_j
=\sum_jg(\e^{g'}_{ij})ge_j\\&=\sum_j\sum_kg(\e^{g'}_{ij})\e^g_{jk}e_k
=\sum_j(\sum_kg(\e^{g'}_{ik})\e^g_{kj})e_j.
\end{align*}
So we can conclude that
\begin{equation*}
\e^{gg'}_{ij}=\sum_kg(\e^{g'}_{ik})\e^g_{kj}
\end{equation*}
\end{proof}

From this result we immediately have

\begin{corollary}
$(\e^g)^{-1}=g(\e^{g^{-1}})$.
\end{corollary}

\subsection{Coordinate description of tensor operations}

We will now investigate how the connections of GF-difference equations behave when we 
perform the usual linear algebra operations on the corresponding modules. 

Let $\e$ and $\fc$ be left A-modules and let $\{e_i\}_{i=1}^{n},\{f_j\}_{j=1}^{m}$ be k-basises for
$\e$ and $\fc$. Let $\{\e^g\}_{g\in G}$ and $\{\fc^g\}_{g\in G}$ be the connection of the
$G$-difference equations $\e$ and $\fc$ with respect to the given basis.
With respect to direct sum we have the following result

\begin{proposition}
$(\e\ds\fc)^g=\e^g\ds\fc^g$.
\end{proposition}
\begin{proof}
A basis for $\e\ds\fc$ is $\{(e_i,0),(0,f_j)\}_{i=1,j=1}^{n,m}$ and
\begin{align*}
g(e_i,0)&=(ge_i,0)=(\sum_j\e^g_{ij}e_j,0)\\
&=\sum_j\e^g_{ij}(e_j,0),\\
g(0,f_i)&=(0,gf_i)=(0,\sum_j\fc^g_{ij}f_j)\\
&=\sum_j\fc^g_{ij}(0,f_j).
\end{align*}
So we have proved that $(\e\ds\fc)^g=\e^g\ds\fc^g$.
\end{proof}

For k-tensorproduct we have a similar simple result.
\begin{proposition}
$(\e\tpk\fc)^g=\e^g\tpk\fc^g$.
\end{proposition}
\begin{proof}
Since $\e$ and $\fc$ are free k-modules it follows that $\{e_i\tpk f_j\}_{i=1,j=1}^{n,m}$ is a k-basis for
 $\e\tpk\fc$. Furthermore we have
\begin{align*}
g(e_i\tpk f_j)&=ge_i\tpk gf_j=(\sum_r\e^g_{ir}e_r)\tpk (\sum_s\fc^g_{js}f_s)\\
&=\sum_r\sum_s\e^g_{ir}\fc^g_{js}e_r\tpk f_s.
\end{align*}
So 
\begin{equation*}
(\e\tpk\fc)^g_{irjs}=\e^g_{ir}\fc^g_{js}.
\end{equation*}
\end{proof}

For $Hom_k(\e,\fc)$ we have the following result.

\begin{proposition}
$Hom_k(\e,\fc)^g=\fc^g\tpk((\e^g)^t)^{-1}$.
\end{proposition}
\begin{proof}
Define elements $\delta_{ij}\in Hom(\e,\fc)$ by
\begin{equation*}
\delta_{ij}(e_k)=\begin{cases}
	f_i& \text{if $j=k$},\\
	0&   \text{if $j\ne k$.}
	\end{cases}
\end{equation*}

Then $\{\delta_{ij}\}$ is a k-basis for $Hom_k(\e,\fc)$. Furthermore we have
\begin{align*}
(g\delta_{ij})(e_u)&=g(\delta_{ij}(g^{-1}e_u))
=g(\delta_{ij}(\sum_v\e^{g^{-1}}_{uv}e_v))\\&=g\sum_v\e^{g^{-1}}_{uv}\delta_{ij}(e_v)
=g(\e^{g^{-1}}_{uj})gf_i\\&=\sum_rg(\e^{g^{-1}}_{uj})\fc^g_{ir}f_r.
\end{align*}
So we can conclude that
\begin{equation*}
Hom_k(\e,\fc)^g=\fc^g\tpk(g(\e^{g^{-1}}))^t.
\end{equation*}
But $g\e^{g^{-1}}=(\e^g)^{-1}$ and the proof is complete.
\end{proof}

Consider the special case $\e^*=Hom_k(\e,k)$. Clearly $\fc^g=k^g=1$ and the previous proposition gives

\begin{corollary}
$(\e^*)^g=((\e^g)^t)^{-1}$.
\end{corollary}

For symmetric and antisymmetric product there are no simple general formulas for computing the connection of
$S^n\e$ and $\wedge^n\e$ in terms of the connection of $\e$. For notational simplicity we only consider the case
$n=2$. Let $\{e_i\}$ be a k-basis for $\e$. Then $\{e_ie_j\}_{i\leqslant j}$ is a k-basis for $S^2\e$. We 
have by definition 
\begin{equation*}
ge_i=\sum_j\e^g_{ij}e_.j
\end{equation*}
But then we have
\begin{align*}
g(e_ie_j)&=ge_ige_j=(\sum_k\e^g_{ik}e_k)(\sum_l\e^g_{jl}e_l)=\sum_{kl}\e^g_{ik}\e^g_{jl}e_ke_l\\
&=\sum_{k>l}\e^g_{ik}\e^g_{jl}e_ke_l+\sum_{k= l}\e^g_{ik}\e^g_{jl}e_ke_l+\sum_{k<l}\e^g_{ik}\e^g_{jl}e_ke_l\\
&=\sum_{k}\e^g_{ik}\e^g_{jk}e_ke_k+\sum_{k<l}(\e^g_{ik}\e^g_{jl}+\e^g_{il}\e^g_{jk})e_ke_l
\end{align*}
So for $n=2$ we have
\begin{equation*}
(S^2\e)^g_{ijkl}=\begin{cases}
	\e^g_{ik}\e^g_{jk}& \text{$k=l$}\\
	\e^g_{ik}\e^g_{jl}+\e^g_{il}\e^g_{jk}& \text{$k<l$}.
	\end{cases}
\end{equation*}

We now consider $\wedge^2\e$. A k-basis for $\wedge^2\e$ is $\{e_i\wedge e_j\}_{i<j}$. Furthermore we have
\begin{align*}
g(e_i\wedge e_j)&=ge_i\wedge ge_j=(\sum_k\e^g_{ik}e_k)\wedge (\sum_l\e^g_{jl}e_l)
=\sum_{kl}\e^g_{ik}\e^g_{jl}e_k\wedge e_l\\
&=\sum_{k<l}\e^g_{ik}\e^g_{jl}e_k\wedge e_l+\sum_{k>l}\e^g_{ik}\e^g_{jl}e_k\wedge e_l\\
&=\sum_{k<l}(\e^g_{ik}\e^g_{jl}-\e^g_{il}\e^g_{jk})e_k\wedge e_l.
\end{align*}
We can conclude that
\begin{equation*}
(\wedge^2\e)^g_{ijkl}=\e^g_{ik}\e^g_{jl}-\e^g_{il}\e^g_{jk}.
\end{equation*}
Assume that $dim_k\e=2$. Then we have
\begin{equation*}
\e^g=\begin{pmatrix} \e^g_{11}&\e^g_{12}\\
		      \e^g_{21}&\e^g_{22}
	\end{pmatrix}
\end{equation*}

The module $\wedge^2\e$ is one dimensional over k with basis $e_1\wedge e_2$. The matrix 
$(\wedge^2\e)^g$ is the scalar
\begin{equation*}
(\wedge^2\e)^g=\e^g_{11}\e^g_{22}-\e^g_{12}\e^g_{21}=det(\e^g).
\end{equation*}
It is evident that a similar result holds in general
\begin{proposition}
Let $dim_k\e=n$. Then
\begin{equation*}
(\wedge^n\e)^g=det(\e^g).
\end{equation*}
\end{proposition}

\subsection{Coordinate description of A-module morphisms and solutions}

Let $\e$ and $\fc$ be GF-difference equations with k-bases $\{e_i\}_{i=1}^n,\{f_i\}_{i=1}^m$
and connection $\{\e^g\}_{g\in G},\{\fc^g\}_{g\in G}$. Let $\phi\in Hom_A(\e,\fc)$ be
a A-module morphism. Then $\phi$ is k-linear and has a matrix $\phi=(\phi_{ij})$ with
respect to the given bases for $\e$ and $\fc$.

\begin{proposition}
$g(\phi)\fc^g=\e^g\phi$ for all $g\in G$.
\end{proposition}
\begin{proof}
The k-morphism $\phi$ is a A-module morphism only if $\phi(ge)=g\phi(e)$ for all $g\in G$.
But we have
\begin{align*}
\phi(ge_i)&=\phi(\sum_j\e^g_{ij}e_j)=\sum_j\e^g_{ij}\phi(e_j)=\sum_j\e^g_{ij}\sum_k\phi_{jk}f_k\\
&=\sum_k(\sum_j\e^g_{ij}\phi_{jk})f_k,\\
g\phi(e_i)&=g(\sum\phi_{ij}f_j)=\sum_jg(\phi_{ij})gf_j=\sum_jg(\phi_{ij})\sum_k\fc^g_{jk}f_k\\
&=\sum_k(\sum_jg(\phi_{ij})\fc^g_{jk})f_k
\end{align*}

comparing sides we have $\sum_jg(\phi_{ij})\fc_{jk}^g=\sum_j\e^g_{ij}\phi_{jk}$ for all i,k
and $g\in G$.
\end{proof}

It is clearly sufficient that the matrix $\phi$ satisfie the equation in the previous 
proposition only on a set of generators for G.
We have defined solution of type $\fc$ to a given GF-difference equation $\e$ as A-morphisms from $\e$ to
$\fc$. Using coordinates as in the previous proposition we see that solutions in our sense is
a matrix of functions on $\set$ that solves a set of classical difference equations. This show
that we are not redefining the notion of solution and justifies our use of this term in our theory.

\section{Invariant structures}

Let $\e,\e',\fc$ and $\fc'$ be GF-difference equations. Let $\phi\in Hom_A(\e,\fc)$ and $\psi\in Hom_A(\e',\fc')$
and define maps $\phi\tpk\psi,S^n\phi,\wedge^n\phi$ by
\begin{align*}
(\phi\tpk\psi)(e\tpk e')&=\phi(e)\tpk\psi(e'),\\
S^n\phi(e_1e_2\cdots e_n)&=\phi(e_1)\phi(e_2)\cdots\phi(e_n),\\
\wedge^n\phi(e_1\wedge e_2\cdots\wedge e_n)&=\phi(e_1)\wedge\phi(e_2)\cdots\wedge(e_n).
\end{align*}

Then the maps are well defined and we have the following result whose proof can
be found in standard texts \cite{lang}.
\begin{proposition}
The maps $\phi\tpk\psi,S^n\phi,\wedge^n\phi$ are well defined and
\begin{align*}
\phi\tpk\psi&\in Hom_A(\e\tpk\e',\fc\tpk\fc'),\\
S^n\phi&\in Hom_A(S^n\e,S^n\fc),\\
\wedge^n\phi&\in Hom_A(\wedge^n\e,\wedge^n\fc).
\end{align*}
\end{proposition}

\subsection{Conserved quantities}

Let $\e ,\fc$ be GF-difference equations and let $\phi\in Hom_A(\e,\fc)$ be
a solution of $\e$ of type $\fc$. Using $\phi$ we can generate morphisms $S^n\phi:S^n\e\ft S^n\fc$ and
$\wedge^n\phi:\wedge^n\phi\ft\wedge^n\fc$. Let us first consider the symmetric case. 
Assume that there exists
an invariant element or structure $\alpha\in S^n\e$. This means that $g\alpha=\alpha$ for all $g\in G$.
Then $S^n\phi(\alpha)\in S^n\fc$ and we have
\begin{align*}
g(S^n\phi(\alpha))=S^n\phi(g\alpha)=S^n\phi(\alpha)
\end{align*}

So that $S^n\phi(\alpha)$ is an invariant structure in $S^n\fc$. Let us consider the particular case when
$\fc$ is the simple object $\s\iso k$ corresponding to the trivial action of G. Then $S^n\s\iso k$ and
we have the  following result
\begin{proposition}
Let $\e$ be any GF-difference equation with an invariant structure $\alpha\in S^n\e$. Let $\phi\in Hom_A(\e,k)$
be any solution of $\e$ of type $\s\iso k$. Then 
\begin{align*}
S^n\phi(\alpha)=constant.
\end{align*}
\end{proposition}
\begin{proof}
We know that $S^n\phi(\alpha)$ is an element in k and that $gS^n\phi(\alpha)=S^n\phi(\alpha)$ for all
elements in G. But G acts transitivly so $S^n\phi(\alpha)$ must be a constant function.
\end{proof}

Note that in coordinates this gives us a symmetric polynomial invariant for the equation $\e$.

Let us next consider the antisymmetric case. Let $\alpha\in\wedge^n\e$ be a invariant structure for $\e$ so
that $g\alpha=\alpha$ for all $g\in G$. Then $g(\wedge^n\phi(\alpha))=\wedge^n\phi(\alpha)$ so $\wedge^n\phi(\alpha)$
is an invariant structure in $\wedge^n\fc$. This leads to the following proposition
 \begin{proposition}
Let $\s$ be a GF-difference equation and $n=dim_k\fc$. Let $\e$ be any 
GF-difference equation with invariant structure $\alpha\in\wedge^n\e$. Then
\begin{align*}
\wedge^n\phi(\alpha)=constant.
\end{align*}
for any solution $\phi\in Hom_A(\e,\s)$.
\end{proposition}

This gives us an antisymmetric polynomial invariant for the equation $\e$. In a similar way conditions
for other types of conservation laws can be specified through invariants.

\subsection{Self-dual equations}

Let $\e$ be any GF-difference equation and let $\e^*$ be the dual equation. Assume that there is
an invariant structure $\alpha\in S^2\e^*$ or $\alpha\in\wedge^2\e^*$. Define a map $F_{\alpha}:\e\ft\e'$ by
\begin{align*}
F_{\alpha}(e)(e')=\alpha(e,e').
\end{align*}
We have the following result
\begin{proposition}
$F_{\alpha}$ is a A-morphism.
\end{proposition}
\begin{proof}
It is evident that $F_{\alpha}(e)\in\e^*$ and that $F_{\alpha}$ is k-linear.
Furthermore we have
\begin{align*}
F_{\alpha}(ge)(e')&=\alpha(ge,e')=gg^{-1}(\alpha(ge,g(g^{-1}e')))=g((g^{-1}\alpha)(e,g^{-1}e'))\\
&=g(\alpha(e,g^{-1}e'))=g(F_{\alpha}(e)(g^{-1}e'))=(gF_{\alpha}(e))(e').
\end{align*}
\end{proof}

Let us now define the notion of self duality for GF-difference equations.
\begin{definition}
A GF-difference equation $\e$ is self-dual if $\e\iso \e^*$ as A-modules.
\end{definition}

Using the previous proposition we can now prove the following
\begin{proposition}
Let $\e$ be a GF-difference equation and assume $\e$ has a nondegenerate  invariant structure $\alpha\in S^2\e^*$
or $\alpha\in\wedge^2\e^*$.
Then $\e$ is self-dual.
\end{proposition}
\begin{proof}
We have a A-morphism $F_{\alpha}:\e\ft\e^*$. This map is bijective if $\alpha$ is nondegenerate
because $F_{\alpha}(e)=\beta$ if and only if $\alpha(e,e')=\beta(e')$ for all $e'$ and these equations
has one and only one solution $e$ since $\alpha$ is nondegenerate.
\end{proof}

This proposition show that any GF-difference equation with an invariant euclidian or symplectic structure
is self-dual.

\subsection{Solutions and composition principles}

Let $\e,\fc$ be GF-difference equations and let $\phi\in Hom_k(\e,\fc)$.
Then $\phi$ is a solution of type $\fc$ of $\e$ if $g\phi=\phi$ for all $g\in G$. This means that
a solution is a invariant structure in $Hom_k(\e,\fc)$. Let $\alpha\in S^2(Hom_k(\e,\fc)^*)\tpk Hom_k(\e,\fc)$
be a invariant structure. Using the standard isomorphism $Hom_k(\fc,\fc')\iso \fc^*\tpk\fc'$,
 $\alpha$ defines a map $T_{\alpha}:Hom_k(\e,\fc)\tpk Hom_k(\e,\fc)\ft Hom_k(\e,\fc)$
defined by
\begin{align*}
T_{\alpha}(\phi,\psi)=\alpha(\phi,\psi)
\end{align*}
For the map $T_{\alpha}$ we have the following result
\begin{proposition}
Let $\e,\fc$ be a GF-difference equations and. Let $\phi,\psi \in Hom_A(\e,\fc)$ be a
pair of solutions of $\e$ of type $\fc$. Then $T_{\alpha}(\phi,\psi)\in Hom_A(\e,\fc)$ is a solution
of $\e$ of type $\fc$.
\end{proposition}
\begin{proof}
We have $g\phi=\phi$ and $g\psi=\psi$ for all $g\in G$ since they are solutions. But then we have
\begin{align*}
g(T_{\alpha}(\phi,\psi))&=g(\alpha(g^{-1}g\phi,g^{-1}g\psi))=(g\alpha)(g\phi,g\psi)\\
&=\alpha(\phi,\psi)=T_{\alpha}(\phi,\psi).
\end{align*}
\end{proof}

So a GF-difference equation $\e$ has a symmetric composition principle for solutions of type $\fc$ 
if there is an invariant structure in $S^2(Hom_k(\e,\fc)^*)\tpk Hom_k(\e,\fc)$. In a similar way
other types of composition principles will correspond to the existense of certain invariants in the
tensor algebra of the equation $\e$.

\section{Module description of classical difference equations}

We will now develope the analog of differential operators on sections in
vectorbundles. Many of the constructions introduced also applies in the case of equations
that are not of finite type. We will however in this section assume that all modules that
appears are GF-difference equations. This will in particular mean that A itself must be
a GF-difference equation. This can only happend if G is a finite group. Since G acts
transitively on $\set$ this means that we are considering the situation where $\set$ is
a finite set.

\subsection{The module of difference operators}
Let $\e$,$\e'$ be A-modules. We will define an action of elements in the module $Hom_k(\e,\e')\tpk A$
on $\e$. For each element $(\phi,a)\in Hom_k(\e,\e')\times A$ define a map $\mu(\phi,a):\e\ft\e'$ by
\begin{align*}
\mu(\phi,a)(e)=\phi(ae).
\end{align*}

\begin{proposition}
$\mu(\phi,a)$ is $\f$-linear.
\end{proposition}
\begin{proof}
\begin{align*}
\mu(\phi,a)(e+e')&=\phi(a(e+e'))=\phi(ae+ae')=\phi(ae)+\phi(ae')\\
&=\mu(\phi,a)(e)+\mu(\phi,a)(e'),\\
\mu(\phi,a)(re)&=\phi(a(re))=\phi(r(ae))=r\phi(ae)=r\mu(\phi,a)(e).
\end{align*}
\end{proof}

Let $\mu$ be the map $(\phi,a)\ft\mu(\phi,a)$. Then we have

\begin{proposition}
$\mu$ is k-bilinear.
\end{proposition}
\begin{proof}
\begin{align*}
\mu(\phi+\phi',a)(e)&=(\phi+\phi')(ae)=\phi(ae)+\phi'(ae)\\
&=\mu(\phi,a)(e)+\mu(\phi',a)(e),\\
\mu(\phi,a+a')(e)&=\phi((a+a')e)=\phi(ae+a'e)=\phi(ae)+\phi(a'e)\\
&=\mu(\phi,a)(e)+\mu(\phi,a')(e),\\
\mu(f\phi,a)(e)&=(f\phi)(ae)=f(\phi(ae))=\phi(f(ae))\\
&=\phi((fa)e)=\mu(\phi,fa)(e).
\end{align*}
\end{proof}
So we have a well defined map $\mu:Hom_k(\e,\e')\tpk A\ft Hom_{\f}(\e,\e')$ defined by
\begin{align*}
\mu(\phi\tpk a)=\phi(ae)
\end{align*}

\begin{proposition}
The map $\mu$ is a A-module morphism.
\end{proposition}
\begin{proof}
By construction the map $\mu$ is a k-module morphism. Let $g\in G$, then we have
\begin{align*}
\mu(g (\phi\tpk a))(e)&=\mu(g\phi\tpk ga)=(g\phi)(gae)\\&=g(\phi(g^{-1}(gae)))=g(\phi(ae))=g(\mu(\phi\tpk a)(e))\\
&=(g\mu(\phi\tpk a))(e).
\end{align*}
\end{proof}

The elements $\theta\in Hom_k(\e,\e')\tpk A$ thus acts as $\f$-linear maps from the module $\e$ to the
module $\e'$. The action is defined by
\begin{align*}
\theta(e)=\mu(\theta)(e).
\end{align*}
 We will now consider the coordinate expression for these maps. Let $\{e_i\},\{e_i'\}$
be k-basises for $\e$ and $\e'$. Then $\{\phi_{ij}\}$ is a basis for $Hom_k(\e,\e')$ where
$\phi_{ij}(e_k)=\delta_{ik}e_j'$.
Let $\theta\in Hom_k(\e,\e')\tpk A$ and $e\in \e$. Using basis we have  $\theta=\sum_{ijg}
\theta_{ijg}\phi_{ij}\tpk g$ and $e=\sum_if_ie_i$. This gives us
\begin{align*}
(\sum_{ijg}\theta_{ijg}\phi_{ij}\tpk g)(\sum_ke_k)&=\sum_{ijg}\theta_{ijg}\phi_{ij}(g\sum_kf_ke_k)
=\sum_{ijg}\theta_{ijg}\phi_{ij}(\sum_kg(f)_k\sum_l\e^g_{kl}e_l)\\
&=\sum_{ijklg}\theta_{ijg}g(f)_k\e^g_{kl}\phi_{ij}(e_l)
=\sum_{ijklg}\theta_{ijg}g(f)_k\e^g_{kl}\delta_{il}e'_j\\
&=\sum_{ijkg}\theta_{ijg}g(f)_k\e^g_{ki}e_j'.
\end{align*}

The equation $\theta(e)=0$ is therefore equivalent to a system of classical difference equations
\begin{align*}
\sum_k(\sum_{ig}\theta_{ijg}\e^g_{ki}g)f_k=0
\end{align*}

In general any elements in the kernel of $\mu$ will be trivial when
considered as $\f$-linear maps.
\begin{example}
Let $\set=\{x,y,z\}$ be the cyclic graph of three elements with symmetry group $S_3$. Let
the group elements in cycle notation be $g_0=id,g_1=(1,3,2),g_2=(1,2,3),g_3=(1,2),g_4=(2,3)$
and $g_5=(1,3)$. Then the element $\theta=\phi\tpk(g_0+g_1+g_2-g_3-g_4-g_5)$, $\phi\ne 0$, is trivial
as a $\f$-linear map.
\end{example}

We therefore makes the following definition
\begin{definition}
$\df(\e,\e')=(Hom_k(\e,\e')\tpk A)/ker\mu$ is the module of difference operators from $\e$ to $\e'$.
\end{definition}

\subsection{Composition of difference operators}
 
Let $\e_1,\e_2$ and $\e_3$ be A-modules. For each pair of elements $(\phi,g)\in Hom_k(\e_2,\e_3)\times A$
define a map $F_g^{\phi}:Hom_k(\e_1,\e_2)\times A\ft Hom_k(\e_1,\e_3)\tpk A$ by
\begin{align*}
F_g^{\phi}(\psi,b)=\phi\circ g\psi\tpk gb
\end{align*}

\begin{proposition}
$F_g^{\phi}$ is middle k-linear for each $(\phi,g)\in Hom_k(\e_2,\e_3)\times A$.
\end{proposition}
\begin{proof}
\begin{align*}
F_g^{\phi}(\psi+\psi',b)&=\phi\circ(g(\psi+\psi'))\tpk gb=\phi\circ(g\psi+g\psi')\tpk gb=(\phi\circ g\psi+
\phi\circ g\psi')\tpk gb\\&=\phi\circ g\psi\tpk gb+\phi\circ g\psi'\tpk gb=F_g^{\phi}(\psi,b)+F_g^{\phi}
(\psi',b),\\
F_g^{\phi}(\psi,b+b')&=\phi\circ g\psi\tpk g(b+b')=\phi\circ g\psi\tpk (gb+gb')\\&=\phi\circ g\psi\tpk gb
+\phi\circ g\psi\tpk gb'=F_g^{\phi}(\psi,b)+F_g^{\phi}(\psi,b'),\\
F_g^{\phi}(\psi,fb)&=\phi\circ g\psi\tpk g(fb)=\phi\circ g\psi\tpk g(f)gb=g(f)(\phi\circ g\psi)\tpk gb\\
&=\phi\circ(g(f)(g\phi))\tpk gb=\phi\circ((g(f)g)\psi)\tpk gb=\phi\circ((g(f))\psi)\tpk gb\\&=\phi\circ g(f\psi)\tpk gb=F_g^{\phi}(f\psi,b).\\ 
\end{align*}
\end{proof}

So we have a well defined $\f$-linear map $F_g^{\phi}:Hom_k(\e_1,\e_2)\tpk A\ft Hom_k(\e_1,\e_3)\tpk A$ defined by
\begin{align*}
F_g^{\phi}(\psi\tpk b)=\phi\circ g\psi\tpk gb.
\end{align*}
We use this map to define a map $F:Hom_k(\e_2,\e_3)\times A\ft Hom_{\f}(Hom_k(\e_1,\e_2)\tpk A,
Hom_k(\e_1,\e_3)\tpk A)$ by
\begin{align*}
F(\phi,a)=\sum_ga_gF_g^{\phi}
\end{align*}
where $a=\sum_ga_gg$.
\begin{proposition}
F is k-bilinear.
\end{proposition}
\begin{proof}
We have
\begin{align*}
F_g^{\phi+\phi'}(\psi\tpk b)&=(\phi+\phi')\circ g\psi\tpk gb=(\phi\circ g\psi+\phi'\circ g\psi)\tpk gb\\
&=\phi\circ g\psi\tpk gb+\phi'\circ g\psi\tpk gb=F_g^{\phi}(\psi\tpk b)+F_g^{\phi'}(\psi\tpk b).
\end{align*}
Using this we find
\begin{align*}
F(\phi+\phi',a)&=\sum_ga_gF_g^{\phi+\phi'}=\sum_ga_gF_g^{\phi}+\sum_ga_gF_g^{\phi'}=F(\phi,a)+F(\phi',a),\\
F(\phi,a+a')&=\sum_g(a_g+a'_g)F_g^{\phi}=\sum_ga_gF_g^{\phi}+\sum_ga'_gF_g^{\phi}=F(\phi,a)+F(\phi,a'),\\
F(f\phi,a)(\psi\tpk b)&=\sum_ga_gF_g^{f\phi}(\psi\tpk b)=\sum_ga_g((f\phi)\circ g\psi\tpk gb)\\
&=\sum_ga_g(f(\phi\circ g\psi)\tpk gb)=\sum_ga_gf(\phi\circ g\psi\tpk gb)\\
&=\sum_g(fa_g)(\phi\circ g\psi\tpk gb).
=F(\phi,fa)(\psi\tpk b)
\end{align*}
\end{proof}

We can conclude that we have a well defined map $F:Hom_k(\e_2,\e_3)\tpk A\ft Hom_{\f}(Hom_k(\e_1,\e_2)\tpk A,
Hom_k(\e_1,\e_3)\tpk A)$ 
\begin{align*}
F(\phi\tpk a)(\psi\tpk b)=\sum_ga_g \phi\circ(g \psi)\tpk gb
\end{align*}
Let $\mu_i,i=1,2,3$ be the action map. Then we have
\begin{proposition}
$F(\theta_2)(\theta_1)\in ker\mu_3$ if $\theta_1\in ker\mu_1$ or $\theta_2\in ker\mu_2$
\end{proposition}
\begin{proof}
Let $\theta_1\in Hom_k(\e_1,\e_2)\tpk A,\theta_2\in Hom_k(\e_2,\e_3)\tpk A$. Then $\theta_1=\sum_i\psi_i\tpk b_i$
and $\theta_2=\sum_i\phi_i\tpk a_i$. But then we have
\begin{align*}
F(\theta_2)(\theta_1)&=\sum_{ij}F(\phi_i\tpk a_i)(\psi_j\tpk b_j)=\sum_{ijg}a_{ig}\phi_i\circ g\psi_j\tpk gb_j
\end{align*}
Using this we find
\begin{align*}
\mu_3(F(\theta_2)(\theta_1))(e_1)&=\sum_{ijg}a_{ig}(\phi\circ g\psi_j)(gb_je_1)=\sum_{ijg}a_{ig}\phi_i
(g\psi_j(b_je_1))\\&=\sum_{ig}a_{ig}\phi_i(g(\sum_j\psi_j(b_je_1)))=\sum_{ig}a_{ig}\phi_i(g\mu_1(\theta_1)(e_1)
\\&=\sum_i(\phi_i(\sum_ga_{ig}\mu_1(\theta_1)(e_1)))=\sum_i\phi_i(a_i\mu_1(\theta_1)(e_1))\\&=\mu_2(\theta_2)(\mu_1(\theta_1)(e_1))=0
\end{align*}
if $\theta_1\in ker\mu_1$ or $\theta_2\in ker\mu_2$
\end{proof}
The map F therefore restricts to the modules of difference operators and we have a map $c:\df(\e_2,\e_3)\times
\df(\e_1,\e_2)\ft\df(\e_1,\e_3)$ defined by
\begin{align*}
c([\theta_2],[\theta_1])=[F(\theta_2)(\theta_1)].
\end{align*}

We use the map c to define composition of difference operators.
\begin{definition}
Let $\Delta_1\in\df(\e_1,\e_2)$ and $\Delta_2\in\df(\e_2,\e_3)$ be difference operators. Define the composition
$\Delta_2\circ\Delta_1\in\df(\e_1,\e_3)$ by
\begin{align*}
\Delta_2\circ\Delta_1=c(\Delta_2,\Delta_1)
\end{align*}
\end{definition}

\subsection{Modules corresponding to difference operators}

By construction $\df(k,k)$ is a left A-module. Let $\mu:A\ft Hom_{\f}(k,k)$ be the action map of A.
Then by definition $\df(k,k)=A/ker\mu$. For any element $a\in A$ let $\Delta_a=[id\tpk a]\in\df(k,k)$
be the corresponding difference operator. Then we have
\begin{proposition}
Let $\Delta\in\df(\e,k)$. Then 
\begin{align*}
a\Delta=\Delta_a\circ\Delta.
\end{align*}
\end{proposition}
\begin{proof}
We only need to consider a generating set for $\df(\e,k)$. Let $\Delta=[\alpha\tpk b]$. Then we have
\begin{align*}
a\Delta&=[a(\alpha\tpk b)]=[\sum_ga_gg(\alpha\tpk b)]=[\sum_ga_g(g\alpha\tpk gb)]\\
&=[\sum_ga_gF_g^{id}(\alpha\tpk b)]=[F(id\tpk a)(\alpha\tpk b)]=c([id\tpk a],[\alpha\tpk b])\\
&=\Delta_a\circ\Delta.
\end{align*}
\end{proof}

Let $\Delta\in\df(\e_1,\e_2)$. Define a map $\phi^{\Delta}:\df(\e_2,k)\ft\df(\e_1,k)$ by
\begin{align*}
\phi^{\Delta}(\nabla)=\nabla\circ\Delta
\end{align*}
Then we have
\begin{proposition}
$\phi^{\Delta}$ is a left A-module morphism.
\end{proposition}
\begin{proof}
\begin{align*}
\phi^{\Delta}(a\nabla)&=(a\nabla)\circ\Delta=(\Delta_a\circ\nabla)\circ\Delta=\Delta_a\circ(\nabla\circ\Delta)
=a\phi^{\Delta}(\nabla)
\end{align*}
\end{proof}
\begin{definition}
Let $\Delta\in\df(\e_1,\e_2)$. The GF-difference equation corresponding to $\Delta$ is
\begin{align*}
\e_{\Delta}=Coker\phi^{\Delta}.
\end{align*}
\end{definition}

\subsection{Classical solutions}

Let $\Delta\in\df(\e_1,\e_2)$ be a difference operator. Define the set of classical solutions
$C(\Delta)$ of $\Delta$ by
\begin{definition}
$C(\Delta)=\{e\in\e_1\mid \Delta(e)=0\}$.
\end{definition}

Let $\s_0\iso k$ be the simple module corresponding to trivial action of G.
For each $e\in C(\Delta)$ define a map $\phi_e:\e_{\Delta}\ft \s_0$ by
\begin{align*}
\phi_e([\lambda])=\lambda(e).
\end{align*}
\begin{proposition}
$\phi_e$ is well defined for each $e\in C(\Delta)$.
\end{proposition}
\begin{proof}
Assume that $[\lambda]=[\lambda']$. Then $\lambda-\lambda'=\phi^{\Delta}(\nabla)$ for some
$\nabla\in\df(\e_2,k)$. But then we have
\begin{align*}
\phi_e([\lambda])&=\lambda(e)=\lambda'(e)+\phi^{\Delta}(\nabla)(e)\\&=\lambda'(e)+\nabla(\phi_{\Delta}(e))
=\lambda'(e).
\end{align*}
\end{proof}
\begin{proposition}
$\phi_e\in Hom_A(\e_{\Delta},\s_0)$
\end{proposition}
\begin{proof}
\begin{align*}
\phi_e(a[\lambda])&=\phi_e([a\lambda])=(\Delta_a\circ\lambda)(e)=\Delta_a(\lambda(e))\\
&=(id\tpk a)(\lambda(e))=a\lambda(e)=a\phi_e([\lambda]).
\end{align*}
\end{proof}

Now define a map $\phi:C(\Delta)\ft Hom_A(\e_{\Delta},\s_0)$ by $\phi(e)=\phi_e$. Then we have
\begin{proposition}
$\phi:C(\Delta)\ft\phi(C(\Delta))\sub Hom_A(\e_{\Delta},\s_0)$ is a isomorphism of $\f$-vectorspaces.
\end{proposition}
\begin{proof}
Assume that $\phi(e)=\phi(e')$. Then we have that $\phi_e([\lambda])=\phi_{e'}([\lambda])$ so that
$\lambda(e-e')=0$ for all $\lambda\in\df(\e_1,k)$. Let $\{e_i\}$ be a basis for $\e_1$ and
$\{e_i^*\}$ the dual basis. Then $e-e'=\sum_if_ie_i$, $\lambda_j=[e_j^*\tpk 1]\in\df(\e_1,k)$ and 
we have $f_j=e_j^*(\sum_if_ie_i)=\lambda_j(e-e')=0$. So $e=e'$ and $\phi$ is injective. Furthermore
we have $\phi(re)([\lambda])=\lambda(re)=r\lambda(e)=r\phi(e)([\lambda])$ so $\phi$ is $\f$-linear.
\end{proof}

The previous proposition show that any classical solution of a difference operator $\Delta$
is contained in the set of solutions of $\e_{\Delta}$ of type $\s_0$.

\subsection{Modules corresponding to systems of difference equations}

Any system of difference equations on the space $\set$  is of the form
\begin{align*}
\sum_{k=1}^n(\sum_gc_{kg}^jg)f_k=0\quad\text{for $j=1\cdots m$}
\end{align*}
The given system of difference equations will only fix the k-module structure of the A-modules
$\e_1$ and $\e_2$. It will not fix the A-module structure  or the operator $\Delta\in\df(\e_1,\e_2)$ separately
but will fix a relation between the A-module structure on $\e_1$ and the operator $\Delta$. The space
of solutions of the given system of difference equations must be equal to $C(\Delta)$ Using bases 
$\{e_i\}$ and $\{f_i\}$ for $\e_1$ and $\e_2$ we have  $\Delta=[\sum_{ijg}\theta_{ijg}\phi_{ij}\tpk g]$
and the relation is
\begin{align*}
\sum_i\theta_{ijg}\e_{ki}^g=c^j_{kg}
\end{align*}
where $\e^g$ is the connection for the action of g on $\e_1$. This means that in general we have many different
modules $\e_{\Delta}$ corresponding to a given system of difference equations. However for all these modules
$\e_{\Delta}$ we have $C(\Delta)\in Hom_A(\e_{\Delta},\s_0)$ so they all contain the set of solutions of
the given system of difference equations.

\end{document}